%
%

\magnification\magstephalf
\baselineskip13pt 
\vsize23.5truecm 
\input nils.sty 

\def\quotationone{\smallrm Hjort and Pollard}
\def\quotationtwo{\smallrm Asymptotics for minimisers}
\def\today{May 1993}

\def\sumin{\sum_{i\le n}}
\def\maxin{\max_{i\le n}}
\def\qt#1{\qquad\hbox{\rm#1}}


\centerline{\bigbf Asymptotics for minimisers of convex processes}
\medskip
\centerline{\bf Nils Lid Hjort and David Pollard}

\smallskip
\centerline{\bf University of Oslo and Yale University} 


\smallskip
{{\smallskip\narrower\noindent\baselineskip12pt 
{\csc Abstract.} 
By means of two simple convexity arguments 
we are able to develop a general method for 
proving consistency and asymptotic normality 
of estimators that are defined by minimisation of convex
criterion functions. This method is then applied to 
a fair range of different statistical estimation problems,
including Cox regression, logistic and Poisson regression,
least absolute deviation regression outside model conditions,  
and pseudo-likelihood estimation for Markov chains.

Our paper has two aims.
The first is to exposit the method itself, 
which in many cases, under reasonable regularity conditions,  
leads to new proofs that are simpler than the traditional proofs.
Our second aim is to exploit the method to its limits
for logistic regression and Cox regression, 
where we seek asymptotic results under
as weak regularity conditions as possible. 
For Cox regression in particular we are able to 
weaken previously published regularity conditions substantially.

\smallskip\noindent
{\csc Key words:}
{\sl argmin lemma approximation, 
convexity, 
Cox regression, 
LAD regression,
log-concavity, 
logistic regression,
minimal conditions, 
partial likelihood, 
pseudo-likelihood}
\smallskip}}

\bigskip
{\bf 1. Introduction.}
This paper develops a simple method for proving consistency
and asymptotic normality for estimators defined by minimisation
of a convex criterion function. Versions of 
the method have been used or partially used 
by several authors, for various specific occasions, 
including 
Jure\v{c}kov\'a~(1977, 1991), 
Andersen and Gill (1982), 
Hjort (1986, 1988a), 
Haberman (1989), Pollard (1990, 1991),
Bickel, Klassen, Ritov and Wellner (1992),
Niemiro (1992), but the general principle has not been widely recognised.

Our aims in this paper are twofold. 
(i) The primary objective is to explain the basic method,
and to illustrate its use in a fair range of statistical
estimation problems. 
In section 2 we state and prove some general theorems about estimators 
that are defined via some form of convex minimisation,    
and in sections 3 and 4 illustrate their use by means of applications to  
sample quantiles, maximum likelihood estimation
when the likelihood is log-concave, and least squares and least
absolute deviation linear regression outside model conditions. 
Similarly sections 5 and 6 treat logistic and Cox regression,
while still further applications are reported in section 7,
including Poisson regression. 
The proofs are relatively simple and instructive,
at least when regularity conditions are kept reasonable. 
(ii) The second objective is to improve on previously 
published results, in the sense of pruning 
down the regularity conditions of theorems 
for two important models, namely logistic regression 
in section 5 and Cox regression in sections 6 and 7A.  
The two aims are mildly conflicting, editorially speaking. 
We soften the conflict in sections 5 and 6 
by writing down first a simple version of a theorem with a simple proof,
and then a harder version with a harder proof.
In this way we hope that
our article has some pedagogic merits while at the same time 
also offering something to the specialists. 

Instead of treating minimisation as a search for a root of a derivative, 
we work directly with the argmin (a minimising value)
of a random function
and are able to approximate it with  the argmin  of a simpler random function. 
In this way we manage to avoid special 
arguments that are often used to prove consistency separately.  
Convexity essentially buys us both consistency 
and asymptotic normality with the same 
dollar,  
and sometimes with cheaper regularity conditions. 

The two convexity lemmas that will be used are as follows.

\eject 

{{\smallskip\sl
{\csc Lemma 1: From pointwise to uniform.}
Suppose $A_n(s)$  
is a sequence of convex random functions defined on an open convex set
$\cal S$ in ${\RR}^p$, 
which converges in probability to some $A(s)$, for each $s$. 
Then $\sup_{s\in K}|A_n(s)-A(s)|$
goes to zero in probability, for each compact
subset $K$ of $\cal S$.
\smallskip}}

{\csc Proof:} 
This is proved in Andersen and Gill (1982, appendix), crediting T.~Brown, 
via `diagonal subsequencing' and an appeal to a corresponding
non-stochastic result (see Rockafellar, 1970, Theorem 10.8). 
For a direct  
proof, see Pollard (1991, section 6). \square 

\smallskip
A convex function is continuous and attains it minimum
on compact sets, but it can be flat at its bottom and have several minima. 
For simplicity we speak about `the argmin' when referring
to any of the possible minimisers. The argmin can be 
selected in a measurable way, as explained in Niemiro (1992, p.~1531),
for example. 

{{\smallskip\sl
{\csc Lemma 2: Nearness of argmins.} 
Suppose $A_n(s)$ is convex as in Lemma 1 and is approximated by $B_n(s)$. 
Let $\alpha_n$ be the argmin of $A_n$, and 
assume that $B_n$ has a unique argmin $\beta_n$.
Then there is a probabilistic bound on how far 
$\alpha_n$ can be from $\beta_n$: for each~$\delta>0$,
$${\rm Pr}\{|\alpha_n-\beta_n|\ge \delta\}
	\le {\rm Pr}\{\Delta_n(\delta)\ge \half h_n(\delta)\}, \eqno(1.1)$$
where 
$$\Delta_n(\delta)=\sup_{|s-\beta_n|\le\delta}|A_n(s)-B_n(s)|
	\quad {\sl and} \quad 
  h_n(\delta)=\inf_{|s-\beta_n|=\delta}B_n(s)-B_n(\beta_n). \eqno(1.2)$$
\smallskip}}

{\csc Proof:}
The lemma as stated has nothing to do with convergence
or indeed with the `$n$' subscript at all, of course,
but is stated in a form useful for later purposes. To prove it, 
let $s$ be an arbitrary point outside the ball around $\beta_n$ 
with radius $\delta$, say $s=\beta_n+lu$ for a unit vector 
$u$, where $l>\delta$. Convexity of $A_n$ implies  
$$(1-\delta/l)\,A_n(\beta_n)+(\delta/l)\,A_n(s)
	\ge A_n(\beta_n+\delta u). $$
Writing for convenience $A_n(s)=B_n(s)+r_n(s)$, we deduce
$$\eqalign{
(\delta/l)\,\{A_n(s)-A_n(\beta_n)\}
	&\ge A_n(\beta_n+\delta u)-A_n(\beta_n) \cr
	&=B_n(\beta_n+\delta u)+r_n(\beta_n+\delta u)
	  -B_n(\beta_n)-r_n(\beta_n) \cr
	&\ge h_n(\delta)-2\Delta_n(\delta). \cr}$$ 	
If $\Delta_n(\delta)<\half h_n(\delta)$, then $A_n(s)>A_n(\beta_n)$ 
for all $s$ outside the $\delta$-ball, which means that 
the minimiser $\alpha_n$ must be inside. This proves (1.1). 

It is worth pointing out that any norm on ${\RR}^p$ 
can be used here, and that no assumptions need to be placed
on the $B_n$ function beside the existence of 
the minimiser $\beta_n$.~\square

\smallskip
The two lemmas will deliver more than mere consistency 
when applied to suitably rescaled and 
recentred versions of convex processes.

\smallskip
We record a couple of useful implications of Lemma 2. 
If $A_n-B_n$ goes to zero uniformly on bounded
sets in probability and $\beta_n$ is stochastically bounded, then 
$\Delta_n(\delta)\arr_p0$ by a simple argument.
It follows that $\alpha_n-\beta_n\arr_p0$ 
provided only that $1/h_n(\delta)$ is stochastically bounded for each
fixed $\delta$. This last requirement says that $B_n$ shouldn't flatten out 
around its minimum as $n$ increases. 

{{\smallskip\sl
{\csc Basic Corollary.} 
Suppose $A_n(s)$  is convex and can be represented as
$\half s'Vs+U_n's+C_n+r_n(s)$, where $V$ is symmetric and positive definite,  
$U_n$ is stochastically bounded, 
$C_n$ is arbitrary, 
and $r_n(s)$ goes to zero in probability for each $s$. 
Then $\alpha_n$, the argmin  of $A_n$, is only $o_p(1)$ away 
>from $\beta_n=-V^{-1}U_n$, the argmin of $\half s'Vs+U_n's+C_n$.
If also $U_n\arr_dU$ then  
$\alpha_n\arr_d-V^{-1}U$. 
\smallskip}}

\eject

{\csc Proof:} 
The function $A_n(s)-U_n's-C_n$ is convex and goes to 
$\half s'Vs$ in probability for each $s$. By the first lemma
the convergence is uniform on bounded sets. Let 
$\Delta_n(\delta)$ {be the supremum of $|r_n(s)|$} over
$\{|s-\beta_n|\le\delta\}$. Then, by Lemma~2,  
$$\alpha_n=-V^{-1}U_n+\eps_n,
	\quad {\rm where\ }
  {\rm Pr}\{|\eps_n|\ge\delta\}
	\le {\rm Pr}\{\Delta_n(\delta)\ge\half k\delta^2\}\to0. \eqno(1.3)$$
Here $k$ is the smallest eigenvalue of $V$, 
and $\Delta_n(\delta)\arr_p0$, by the arguments used above.~\square

\smallskip
A useful slight extension of this is when 
$A_n(s)=\half s'V_ns+U_n's+C_n+r_n(s)$ is convex, with a nonnegative definite 
symmetric $V_n$ matrix that converges in probability to 
a positive definite $V$. Writing $V_n=V+\eta_n$ the remainder 
$\eta_n$ can be absorbed into $r_n(s)$ and the result above holds.

\bigskip
{\bf 2. General results for convex minimisation estimators.}
This section presents three basic theorems about the 
asymptotic behaviour of estimators that are defined 
by minimisation of some convex criterion function. 
The first is for the independent identically distributed (i.i.d.)~case. 
The second is stated for independent observations with different
distributions, and is suitable for proving consistency and 
asymptotic normality in regression models, for example,
under model conditions. 
The third theorem also applies to regression model estimators, 
but is suited to give asymptotic results also outside model conditions. 
Applications and illustrations are provided in sections 3, 4 and 5. 


\subsection
{\sl 2A. A theorem for the i.i.d.~case.}
Let $Y_1,Y_2,\ldots$ be i.i.d.~from some distribution $F$. 
A certain $p$-dimensional parameter $\theta_0=\theta(F)$ is of interest. 
Assume that one of the ways of characterising 
this parameter is to say that it minimises 
$Eg(Y,t)=\int g(y,t)\,\d{}F(y)$, where the $g(y,t)$ function is convex in $t$. 
Examples include quantiles, the mean, M-estimation and 
maximum likelihood estimation parameters and so on; see sections 3 and 4.  
In the expectation expression above, and later on, 
$Y$ denotes a generic observation from the true underlying $F$. 

Some weak expansion of $g(y,t)$ around the value $\theta_0$ 
of $t$ is needed, but we avoid explicitly requiring pointwise derivatives
to exist. With this in mind, write 
$$g(y,\theta_0+t)-g(y,\theta_0)=D(y)'t+R(y,t) \eqno(2.1) $$ 
for a $D(y)$ with mean zero under $F$.
If $\E R(Y,t)^2$ is of order $o(|t|^2)$ as $t\arr0$, 
as we will usually require, 
then $D(y)$ is nothing but the derivative in quadratic mean 
of the function $g(y,\theta_0+t)$ at $t=0$. 
  
{{\smallskip\sl
{\csc Theorem 2.1.} 
Suppose that $g(y,t)$ is convex in $t$ as above, and that (2.1) holds with 
$$\E\{g(Y,\theta_0+t)-g(Y,\theta_0)\}
	=\E R(Y,t)
	=\half t'Jt+o(|t|^2) {\sl\ as\ }t\arr0 \eqno(2.2) $$
for a positive definite matrix $J$. 
Suppose also that ${\rm Var}\,R(Y,t)=o(|t|^2)$, and 
that $D(Y)$ has a finite covariance matrix $K=\int D(y)D(y)'\,\d{}F(y)$.
Then the estimator $\hatt\theta_n$ which minimises 
$\sumin g(Y_i,t)$ is $\sqrt{n}$-consistent for $\theta_0$, and 
$$\sqrt{n}(\hatt\theta_n-\theta_0)=
	-J^{-1}n^{-1/2}\sumin D(Y_i)
	+o_p(1). \eqno(2.3)$$
In particular 
$\sqrt{n}(\hatt\theta_n-\theta_0)
\arr_d -J^{-1}\normal_p\{0,K\}=\normal_p\{0,J^{-1}KJ^{-1}\}$. 
\smallskip}}


{\csc Proof:} 
Consider the convex function 
$A_n(s)=\sumin \{g(Y_i,\theta_0+s/\sqrt{n})-g(Y_i,\theta_0)\}$. 
It is minimised by $\sqrt{n}(\hatt\theta_n-\theta_0)$. 
Note first that 
$n\E R(Y,s/\sqrt{n})=\half s'Js+r_{n,0}(s)$ 
where $r_{n,0}(s)=no(|s|^2/n)\arr0$ for fixed $s$. 
Accordingly, using (2.1),  
$$\eqalign{A_n(s)
&=\sumin \bigl\{D(Y_i)'s/\sqrt{n}
	+R(Y_i,s/\sqrt{n})-\E R(Y_i,s/\sqrt{n})\bigr\} 
	+n\E R(Y,s/\sqrt{n}) \cr
&=U_n's+\half s'Js+r_{n,0}(s)+r_n(s), \cr}$$
in which 
$$U_n=n^{-1/2}\sumin D(Y_i)
	\quad {\rm and} \quad 
 r_n(s)=\sumin \{R(Y_i,s/\sqrt{n})-\E R(Y_i,s/\sqrt{n})\}. $$
Now $r_n(s)$ tends to zero in probability for each $s$, 
since its mean is zero and its variance is  
$\sumin {\rm Var}\,R(Y_i,s/\sqrt{n})=no(1/n)$. 
This, together with the Basic Corollary of section 1, 
proves (2.3) and the limit distribution result, 
since $U_n$ goes to a $\normal_p\{0,K\}$ by the central limit theorem. 
Note that both consistency and asymptotic normality followed 
>from the same approximation argument.~\square 


\smallskip
Note that $\Var\,R(Y,t) = \E R(Y,t)^2 + O(t^4)$,  
so we might as well work with second moments rather than variances.
Notice also that the differentiability assumption~(2.2) 
is applied to the process obtained by averaging out over 
the distribution $F$, a smoothing that can
eliminate troublesome pointwise behaviour of $R(y,t)$. Huber~(1967) recognised
this advantage of smoothing before differentiating.

\subsection
{\sl 2B. A theorem for independent observations with
different distributions.} 
Assume that the true density of $Y_i$ is of the form 
$f_i(y_i)=f_i(y_i,\theta_0,\eta_i)$, 
where $\theta_0$ is a certain $p$-dimensional parameter of interest.
Suppose that an estimator $\hatt\theta_n$ for $\theta_0$ is proposed which 
minimises $\sumin g_i(Y_i,\theta)$, where the 
$g_i(y_i,\theta)$ functions are convex in $\theta$. 
A simple example is linear regression, where 
$Y_i=\theta_0'x_i+\eps_i$ and 
$\hatt\theta_n$ minimises $\sumin (Y_i-\theta'x_i)^2$.

Suppose that 
$g_i(y_i,\theta_0+t)-g_i(y_i,\theta_0)=D_i(y_i)'t+R_i(y_i,t)$, 
where $\E D_i(Y_i)=0$. With the previous development in mind, write 
$$\E R_i(Y_i,t)=\half t'A_it+v_{i,0}(t)
	\quad {\rm and} \quad
  {\rm Var}\,R_i(Y_i,t)=v_i(t), \eqno(2.4)$$
and let $B_i$ be the variance matrix for $D_i(Y_i)$.
The sums $J_n=\sumin A_i$ and $K_n=\sumin B_i$ are featured below. 
The first useful result, properly generalising Theorem 2.1, 
is the following,
which is proved by copying the arguments of 2A mutatis mutandis. 

{{\smallskip\sl 
{\csc Theorem 2.2.}
Assume that $\sumin v_{i,0}(s/\sqrt{n})\arr0$
and $\sumin v_i(s/\sqrt{n})\arr0$ for each $s$, 
and that $J_n/n$ and $K_n/n$ converge to $J$ and $K$,
where $J$ is positive definite. Then  
$\sqrt{n}(\hatt\theta_n-\theta_0)$ is only $o_p(1)$
away from $-J^{-1}n^{-1/2}\sumin D_i(Y_i)$. 
If in particular the Lindeberg requirements are fulfilled
for the $D_i(Y_i)$ sequence, then 
$\sqrt{n}(\hatt\theta_n-\theta_0)\arr_d \normal_p\{0,J^{-1}KJ^{-1}\}$.
\smallskip}}

Another result which sometimes is stronger is as follows.
Assume that $\sumin v_{i,0}(J_n^{-1/2}s)\arr0$ and
$\sumin v_i(J_n^{-1/2}s)\arr0$ for each $s$, 
and that $J_n^{-1}K_n$ is bounded. Then 
$$J_n^{1/2}(\hatt\theta_n-\theta_0)
	=-J_n^{-1/2}K_n^{1/2}U_n+o_p(1), \eqno(2.5)$$
where $U_n=K_n^{-1/2}\sumin D_i(Y_i)$. If in particular
there are matrices $J$ and $K$ such that $J_n^{-1}K_n$ goes to $J^{-1}K$,
and the Lindeberg conditions are fulfilled, securing 
$U_n\arr_d\normal_p\{0,I_p\}$, then 
$J_n^{1/2}(\hatt\theta_n-\theta_0)\arr_d\normal_p\{0,J^{-1/2}KJ^{-1/2}\}$. 
This result is proved by studying the convex function 
$\sumin\{g_i(Y_i,\allowbreak\theta_0+J_n^{-1/2}s)-g_i(Y_i,\theta_0)\}$. 
In some situations of interest $J_n=K_n$, further simplifying
(2.5). See section 5 for an illustration of this.

%

\subsection
{\sl 2C. A theorem for regression type estimators 
outside model conditions.} 
The results of 2B are sometimes not sufficient. Theorem 2.3 below will work
for asymptotic behaviour of regression methods outside model
conditions, as made clear in section 3D, for example.

Assume that some covariate vector $x_i=(x_{i,1},\ldots,x_{i,p})'$ 
is associated with observation $Y_i$. 
For simplicity we formulate a result in terms of densities,
rather than general distribution functions. 
Suppose that the true density for $Y_i$ given $x_i$ is $f(y_i|x_i)$
but that some regression model postulates $f(y_i,\beta|x_i)$,
for a suitable $p$-dimensional parameter vector $\beta$. 
We consider an estimator $\hatt\beta_n$ defined to minimise
$\sumin g_i(Y_i,\beta|x_i)$, where $g_i(y_i,\beta|x_i)$ is convex
in $\beta$ for each $(y_i,x_i)$. In the following we shall assume
that the empirical distribution of $x_1,\ldots,x_n$, whether
actually random or under the experimenter's control,
converges to a well-defined distribution $H$ in $x$-space.
This conceptual limit is to be thought of 
as the `covariate distribution'. 
Assume that $n^{-1}\sumin g_i(Y_i,\beta|X_i)$ converges in probability
to a function with a unique minimiser $\beta_0$. 

Under these circumstances it is not generally possible to 
get a representation like the one that led to (2.4), 
because of heterogeneity as well as potential modelling bias, 
as the applications in section 3D and section 5C will illustrate. 
It becomes necessary to include a $x_i$-dependent bias term. 
Suppose that it is possible to write 
$$g_i(y_i,\beta_0+t|x_i)-g_i(y_i,\beta_0|x_i)
	=\{\delta(x_i)+D_i(y_i|x_i)\}'t+R_i(y_i,t|x_i), \eqno(2.6)$$
where $\E D_i(Y_i|x_i)=0$ and ${\rm VAR}\,D_i(Y_i|x_i)=B_i(x_i)$. 
Write furthermore 
$$\E R_i(Y_i,t|x_i)=\half t'A_i(x_i)t+v_{i,0}(t|x_i)
	\quad {\rm and} \quad
  {\rm Var}\,R_i(Y_i,t|x_i)=v_i(t|x_i). \eqno(2.7)$$
This time three matrix sums are needed, 
$J_n=\sumin A_i(x_i)$, $K_n=\sumin B_i(x_i)$, and 
$L_n=\sumin \delta(x_i)\delta(x_i)'$. 

{{\smallskip\sl
{\csc Theorem 2.3.} 
Assume that the $x_1,x_2,\ldots$ sequence is such that 
$$\sumin v_{i,0}(s/\sqrt{n}|x_i)\arr0
	\quad {\sl and} \quad 
  \sumin v_{i}(s/\sqrt{n}|x_i)\arr_p0 \quad {\sl for\ each\ }s, \eqno(2.8)$$
that the $J_n/n$ sequence is bounded away from zero, 
and that the $K_n/n$ and $L_n/n$ sequences are bounded. 
Then 
$$\sqrt{n}(\hatt\beta_n-\beta_0)
	=-(J_n/n)^{-1}\bigl\{n^{-1/2}\sumin \delta(x_i)
	+n^{-1/2}\sumin D_i(Y_i|x_i)\bigr\}+\eps_n, \eqno(2.9)$$
where $\eps_n=\eps_n(x_1,\ldots,x_n)\arr_p0$. 
\smallskip}}

The proof is quite similar to previous proofs in this section,
taking as its starting point the convex function 
$\sumin \{g_i(Y_i,\beta_0+s/\sqrt{n}|x_i)-g_i(Y_i,\beta_0|x_i)\}$.
We omit the details.

\smallskip
The (2.9) representation has two statistically interesting implications.
(i) In the conditional framework with a given $x_i$ sequence, 
suppose that $J_n/n\arr J$ and $K_n/n\arr K$ and that 
the Lindeberg condition holds for $\sumin n^{-1/2}D_i(Y_i|x_i)$. Then 
$$\sqrt{n}(\hatt\beta_n-\beta_0)|x_1,\ldots,x_n 
	=\normal_p \bigl\{-(J_n/n)^{-1}n^{-1/2}\sumin \delta(x_i),
		J^{-1}KJ^{-1}\bigr\}+\eps_n', \eqno(2.10)$$
where $\eps_n'\arr_p 0$. So $\hatt\beta_n$ is approximately 
normal with variance matrix $J^{-1}KJ^{-1}/n$, but 
actually biased with a bias depending on $x_1,\ldots,x_n$.  
The bias is typically zero under exact 
regression model conditions, see 3D below.  
(ii) Secondly, if the $x_i$'s can be treated as being independent and 
coming from their own `design distribution' $H(\d x)$ in $x$-space, 
then $\delta(x_i)$ has mean zero and variance matrix $L$, say. 
In this unconditional framework 
$$\sqrt{n}(\hatt\beta_n-\beta_0)
	=J^{-1}\normal_p \{0,K+L\}+o_p(1)
	\arr_d \normal_p\{0,J^{-1}(K+L)J^{-1}\}. \eqno(2.11)$$  

\bigskip
{\bf 3. Applications and illustrations.}

\subsection
{\sl 3A. The median.}
Let $Y_1,Y_2,\ldots$ be i.i.d.~from a density $f$, 
let $\mu$ be the population median, 
and let $M_n$ be the sample median from the 
first $n$ observations. We shall prove the well known fact  
that $M_n$ is consistent for $\mu$ and that 
$$\sqrt{n}(M_n-\mu)\arr_d\normal\{0,1/4f(\mu)^2\}, \eqno(3.1)$$
provided only that $f$ is positive and continuous at $\mu$. 

This fits into the framework of 2A with 
the convex function $g(y,t)=|y-t|$. 
The (2.1) expansion reads 
$$|y-(\mu+t)|-|y-\mu|=D(y)t+R(y,t), $$
where $D(y)=-I\{y>\mu\}+I\{y\le\mu\}$, and 
$$R(y,t)=\cases{
	2(t-(y-\mu))\,I\{\mu\le y\le \mu+t\} &if $t>0$, \cr
	2((y-\mu)-t)\,I\{\mu+t\le y\le \mu\} &if $t<0$, \cr}$$
while $R(y,0)=0$, which makes it easy to verify 
$$\E R(Y,t)=f(\mu)t^2+o(t^2) 
	\quad {\rm and} \quad 
  \E R(Y,t)^2=\hbox{$4\over3$}f(\mu)|t|^3+o(|t|^3). $$

Actually we only need a distribution function with a positive derivative
at~$\mu$. Of course we don't get the explicit $|t|^3$ bound then.  
Notice that $D(Y)$ and $R(Y,t)$ are bounded functions 
even if $|Y-t|$ itself can have infinite expected value, 
since we work with the difference $|Y-(\mu+t)|-|Y-\mu|$.
Since the variance of $D(Y)$ is equal to 1, assertion (3.1) follows from
Theorem 2.1.  
See 4A below for an extension of this result. 

\subsection
{\sl 3B. Simultaneous asymptotic normality of order statistics.}
Let $f$ be positive and continuous in its support region, and
consider the function 
$$g_p(y,t)=p\{(y-t)_+-y_+\} +(1-p)\{(t-y)_+-(-y)_+\}. $$
It is convex in $t$ and 
its expected value is minimal for $t=F^{-1}(p)=\mu_p$, 
the $p$-th quantile of the underlying distribution, and 
$$\E\{g_p(Y,t)-g_p(Y,\mu_p)\}
	=\half f(\mu_p)(t-\mu_p)^2+o((t-\mu_p)^2) $$
can be shown. The (2.1) expansion works with 
$$D(y)=(1-p)I\{y\le\mu_p\}-pI\{y>\mu_p\}=I\{y\le \mu_p\}-p $$ 
and $\E R(Y,t)^2=O(|t|^3)$ can be checked. 
Let $Q_{n,p}$ be the minimiser of $\sumin g_p(Y_i,t)$,
which is sometimes non-unique, but which in any case is 
at most $O_p(n^{-1})$ away from the $[np]$'th order statistic $Y_{([np])}$. 
The general theorem of 2A implies 
$$Z_n(p)=\sqrt{n}(Q_{n,p}-\mu_p)
	=-f(\mu_p)^{-1}\sqrt{n}\{F_n(\mu_p)-p\}+\eps_n(p), \eqno(3.2)$$
where $F_n$ is the empirical distribution function and 
$\eps_n(p)\arr0$ in probability for each $p$. 
This links the quantile process $Z_n$ to the empirical process, 
and proves finite-dimensional convergence in distribution 
of the quantile process 
to a Gau\ss ian process $Z(.)$ with mean zero and covariance structure 
$${\rm cov}\{Z(p_1),Z(p_2)\}
	={p_1(1-p_2)\over f(\mu_{p_1})f(\mu_{p_2})}
	\quad {\rm for\ }p_1\le p_2. \eqno(3.3)$$
The traditional proofs of this finite-dimensional convergence result 
are rather messier 
than the above. There is in reality also process convergence here, 
of course, which is linked to the fact that 
$\sup_{\delta\le s\le 1-\delta}|\eps_n(p)|$ goes to zero 
in probability for each $\delta$. Proving this is not within 
easy reach of our method, however. See also the comment ending 3D below. 



\subsection
{\sl 3C. Estimation in $L_\alpha$ mode.}
Let more generally $M_{n,\alpha}$ minimise
$\sumin |Y_i-t|^\alpha$, where $\alpha\ge1$,
and let $\xi_\alpha$ be the population parameter
that minimises $\E|Y-t|^\alpha$. For $\alpha={3\over2}$
we would expect an estimator with properties somehow 
between those for the median and the mean, for example. 
We can prove
$$\sqrt{n}(M_{n,\alpha}-\xi_\alpha)
	\arr_d\normal\{0,\tau^2\}
	\quad {\rm where\ }
	\tau^2={\E|Y-\xi_\alpha|^{2(\alpha-1)}
	 \over
	 \{(\alpha-1)\,\E|Y-\xi_\alpha|^{\alpha-2}\}^2}, \eqno(3.4)$$
assuming $\E|Y|^{2(\alpha-1)}$ to be finite.
The proof proceeds by mimicking that for the simpler case $\alpha=1$.
One needs to use 
$$D(y)=-\alpha (y-\xi_\alpha)^{\alpha-1}I\{y>\xi_\alpha\}
	+\alpha(\xi_\alpha-y)^{\alpha-1}I\{y<\xi_\alpha\}, $$
and it is somewhat more cumbersome but feasible 
to bound $\E R(Y,t)^2$.  
And finally needed is the analytical fact that 
$\E\{|Y-(\xi_\alpha+t)|^\alpha
	-|Y-\xi_\alpha|^\alpha\}=\half K_ft^2+o(t^2)$,
in which $K_f=\alpha(\alpha-1)\E|Y-\xi_\alpha|^{\alpha-2}$. 

It is interesting to note here that 
$(\alpha-1)\E|Y-\xi_\alpha|^{\alpha-2}$ 
tends to $2f(F^{-1}(\half))$ as $\alpha$ tends to 1,
explaining the connection from the moment-type expression 
for the variance 
$\tau^2$ of (3.4) to the rather different-looking expression  
for the median case. 

It is also worth pointing out that the (3.4) result can be reached 
via influence functions and function space methods as well. 
The influence function can be found to be 
$$I(F,y)=\cases{
	\alpha K_F^{-1}|y-\xi_\alpha(F)|^{\alpha-1} 
		&if $y>\xi_\alpha(F)$, \cr
	-\alpha K_F^{-1}|y-\xi_\alpha(F)|^{\alpha-1} 
		&if $y<\xi_\alpha(F)$, \cr}$$
after which the usual argument is that since
$\sqrt{n}(M_{n,\alpha}-\xi_\alpha)=n^{-1/2}\sumin I(F,Y_i)+\eps_n,$
for suitable remainder term $\eps_n$, one must have limiting 
normality with $\tau^2=\int I(F,y)^2\,\d{}F(y)$, agreeing with (3.4). 
But proving that $\eps_n$ here goes to zero in probability 
is not trivial, since the $\xi_\alpha$ functional is rather
non-smooth. The argument can be saved via establishing 
Lipschitz differentiability, as in Example 1 of Huber (1967).
Our method manages to avoid these somewhat sophisticated arguments. 


\subsection
{\sl 3D. Agnostic least squares and least absolute deviation regression.} 
Statistical regression is about estimating the unknown 
centre value of $Y$ for given $x$, i.e.~the curve or surface 
${\rm centre}(Y|x)$, based on $p+1$-tuplets $(x_i,Y_i)$,
where `centre' could be the mean or the median. 
Ordinary linear regression uses a linear approximation 
$\beta'x=\sum_{j=1}^p\beta_jx_j$ for this centre function,
which is often a very reasonable method even if the true 
underlying centre function is somewhat non-linear. 
The least squares regression estimator is $\hatt\beta_n'x$ where 
$\hatt\beta_n$ minimises $\sumin (Y_i-\beta'x_i)^2$, and 
the least absolute deviation estimator is $\tilda\beta_n'x$ where 
$\tilda\beta_n$ minimises $\sumin |Y_i-\beta'x_i|$. 

Statistical properties of these estimators are usually
investigated only under the admittedly unlikely assumption that 
the true surface {\it is} linear and that the variances are 
constant over the full region, i.e.
$$Y_i=\beta_0'x_i+\sigma \eps_i \eqno(3.5)$$ 
where the $\eps_i$'s are i.i.d.~standardised residuals
centred around zero. An in some sense more honest approach would 
be to merely postulate that 
$$Y_i=m(x_i)+\sigma(x_i)\eps_i, \eqno(3.6)$$ 
for some smooth functions $m(x)$ and $\sigma(x)$,
and view the regression surface estimator as an 
attempt to produce a good linear approximation to 
the evasive $m(x)$.  
Our plan now is to derive properties under robust and agnostic 
(3.6) conditions using Theorem 2.3 of 2C,
while assuming that the empirical distribution of 
$x_i$'s converges to an appropriate `covariate distribution' $H$. 
Under ideal (3.5) conditions they specialise to 
results obtainable using the simpler Theorem 2.2 of 2B. 

Consider least squares regression first, 
assuming the $\eps_i$'s to have mean zero and variance one. 
This fits into the 2C framework with 
$g_i(Y_i,\beta|x_i)=\half(Y_i-\beta'x_i)^2$. 
The method aims at getting the best linear approximation
$\beta_0'x$ to $m(x)$, in the sense of minimising 
the limit of $n^{-1}\sumin (m(x_i)-\beta'x_i)^2$. 
In fact this means 
$\beta_0=(EXX')^{-1}EXY$. We find 
$$\eqalign{
g_i(Y_i,\beta_0+t|x_i)-g_i(Y_i,\beta_0|x_i)
	&=-(Y_i-\beta_0'x_i)x_i't+\half (t'x_i)^2 \cr
	&=-\bigl(\delta(x_i)+D_i(Y_i|x_i)\bigr)'t+\half t'x_ix_i't, \cr}$$
in which 
$$\delta(x_i)=(m(x_i)-\beta_0'x_i)x_i
	\quad {\rm and} \quad 
  D_i(Y_i|x_i)=(Y_i-m(x_i))x_i. $$ 
In the notation of (2.7) one has $A_i(x_i)=x_ix_i'$ and 
both remainder terms are simply equal to zero. Consider 
$$J_n=\sumin x_ix_i', \quad 
  K_n=\sumin \sigma(x_i)^2x_ix_i', \quad 
  L_n=\sumin \{m(x_i)-\beta_0'x_i\}^2x_ix_i'. \eqno(3.7)$$
Two results can be given, corresponding to (2.10) and (2.11).
First, suppose the $x_i$ sequence is such that 
$J_n/n\arr$ a positive definite $J$, $K_n/n\arr K$, 
that the $L_n/n$ sequence is bounded, and that 
$\maxin \sigma(x_i)^2|x_i|^2/\sumin \sigma(x_i)^2|x_i|^2\arr0$. 
Then $\sqrt{n}(\hatt\beta_n-\beta_0)$ is asymptotically 
normal with mean $J^{-1}n^{-1/2}\sumin (m(x_i)-\beta_0'x_i)x_i$ 
and variance matrix $J^{-1}KJ^{-1}$. 
Secondly, under the unconditional viewpoint where the 
$x_i$'s are seen as i.i.d.~with finite variance 
matrix $L=E(m(X)-\beta_0'X)^2XX'$ for $\delta(x_i)$, then 
$$\sqrt{n}(\hatt\beta_n-\beta_0)\arr_d 
	\normal_p \{0,J^{-1}(K+L)J^{-1}\}. \eqno(3.8)$$
Note that $K+L$ can be estimated consistently with 
$n^{-1}\sumin(Y_i-\hatt\beta_n'x_i)^2x_ix_i'$. 

These results can also be derived more or less directly,
i.e.~without the convex machinery of section 2, see 
Exercise 45 in Hjort (1988b). In the least absolute deviation 
case to be reported on next a direct approach is much more difficult, however,
but it can be efficiently handled using the methods of section 2. 

For the LAD regression case, take 
the $\eps_i$'s of (3.6) to have distribution $F$ 
with median zero and variance one. We will assume that $F$ has a
density $f$ which further possesses a continuous 
derivative $f'$. 
In this case $g_i(Y_i,\beta|x_i)=|Y_i-\beta'x_i|$, and 
the method aims at getting the best approximation
$\beta_0'x$ to $m(x)$ in the sense of minimising 
the long term value of $n^{-1}\sumin \E|m(x_i)-\beta'x_i+\sigma(x_i)\eps_i|$. 
We skip the various details that have to be worked through
to reach a result here. They resemble those above and arguments used in 3A. 
To give the result, let $D_i(Y_i|x_i)=1$ if $Y_i\le\beta_0'x_i$ 
and $-1$ if $Y_i>\beta_0'x_i$, with conditional mean
$h(x_i)=2\,{\rm Pr}\{m(x_i)+\sigma(x_i)\eps_i\le\beta_0'x_i\}-1$, 
and consider the three matrices 
$$J_n=\sumin 2f_i(\beta_0'x_i-m(x_i))x_ix_i', 
	\quad 
  K_n=\sumin \{1-h(x_i)^2\}\,x_ix_i', 
	\quad 
  L_n=\sumin h(x_i)^2\,x_ix_i', $$
where $f_i(z)=f(z/\sigma(x_i))/\sigma(x_i)$ is the density 
of the scaled residual $\sigma(x_i)\eps_i$. 
In particular $K_n+L_n=\sumin x_ix_i'$. 
As for the least squares case these efforts lead to a representation
$$\sqrt{n}(\tilda\beta_n-\beta_0)
	=-(J_n/n)^{-1}\Bigl[n^{-1/2}\sumin h(x_i)x_i
 	 +n^{-1/2}\sumin \{D_i(Y_i|x_i)-h(x_i)\}x_i\Bigr]+\eps_n. \eqno(3.9)$$
This has one implication for given $x_i$-sequences and another
implication for the `overall variability'. Under some 
mild assumptions $J_n/n\arr J$ and $(K_n+L_n)/n\arr K+L$, 
and $\sqrt{n}(\tilda\beta_n-\beta_0)\arr_d \normal_p\{0,
J^{-1}(K+L)J^{-1}\}$. 
The $K+L$ matrix is estimated consistently 
using $\sumin x_ix_i'/n$ whereas a more complicated consistent estimate,
involving smoothing and density estimtaion, 
can be constructed for $J$. 

The special case of ${\rm med}(Y|x)=m(x)=\beta_0'x$ has 
$J_n=\sumin 2f_i(0)x_ix_i'/\sigma(x_i)$, 
and the perfect but perhaps unrealistic case of both
a linear median and a constant variance has 
$J_n^{-1}(K_n+L_n)J_n^{-1}=\{4f(0)^2\}^{-1}(\sumin x_ix_i')^{-1}\sigma^2$. 
This is the case considered in Pollard (1990). 

Our method can also be applied to the quantile regression
situation, where one aims to estimate $m(x_0)+\sigma(x_0)F^{-1}(p)$,
for example, to construct a prediction interval for a future
$Y$ at a given covariate value $x_0$. This time one minimises
$\sumin g_p(Y_i,\beta'x_i)$ with the $g_p$ function of 3B. 
This gives a suitable generalisation of results reached by
Bassett and Koenker (1982). 

\bigskip
{\bf 4. Maximum likelihood type estimation.} 

\subsection
{\sl 4A. Log-concave densities.}
Suppose $Y_1,Y_2,\ldots$ are i.i.d.~from some continuous
density $f$, and that a parametric model of the form 
$f(y,\theta)=f(y,\theta_1,\ldots,\theta_p)$
is employed, where the parameter space is some open and convex region. 
We stipulate that $\log f(y,\theta)$ be concave in $\theta$ in this region
and shall be able to reprove familiar results on
maximum likelihood (ML) and Bayes estimation, 
using the convexity based results of section 2,
but with milder smoothness
assumptions than those traditionally employed. 

Note that the log-likelihood 
$\sumin \log f(Y_i,\theta)$ 
when divided by $n$ tends to 
$\E\log f(Y,\theta)\allowbreak=\int f(y)\log f(y,\theta)\,\d{}y$,
for each $\theta$. Assume that this function has a unique global
maximum at $\theta_0$, which is the `agnostic parameter value'
that gives best approximation according to the Kullback--Leibler
distance $\int f(y)\log\{f(y)/f(y,\theta)\}\,\d{}y$ 
>from truth to approximating density. From 
section 2A the following result is quite immediate.

{{\smallskip\sl
{\csc Theorem 4.1.} 
Suppose 
$\log f(y,\theta_0+t)-\log f(y,\theta_0)=D(y)'t+R(y,t)$
is concave in $t$, for a $D(.)$ function with mean zero and finite
covariance matrix $K$ under $f$, and that the remainder term satisfies
$$\E\{\log f(Y,\theta_0+t)-\log f(Y,\theta_0)\}
      =\E R(Y,t)=-\half t'Jt+o(|t|^2) \eqno(4.1)$$
as well as ${\rm Var}\,R(Y_i,t)=o(|t|^2)$, where $J$ is symmetric
and positive definite. Then the maximum likelihood estimator 
$\hatt\theta_n$ is $\sqrt{n}$-consistent for $\theta_0$ and 
$$\sqrt{n}(\hatt\theta_n-\theta_0)
	=J^{-1}n^{-1/2}\sumin D(Y_i)+o_p(1)
	\arr_d J^{-1}\normal\{0,K\}
	=\normal_p\{0,J^{-1}KJ^{-1}\}. $$
\smallskip}}


In ordinary smooth cases one can Taylor expand
and use $D(y)=\dell\log f(y,\theta_0)/\dell\theta$
and find a remainder $R(y,t)$ with 
mean $-\half t'Jt+O(|t|^3)$ and squared mean of order $O(|t|^4)$,
involving 
$$J=-\E_f{\dell^2\log f(Y_i,\theta_0)\over \dell\theta\dell\theta}
	\quad {\rm and} \quad 
  K={\rm VAR}_f{\dell\log f(Y_i,\theta_0)\over \dell\theta}. \eqno(4.2)$$
Notice that when the model happens to be perfect, 
as in textbooks for optimistic statisticians, 
then $K=J$, and we get the
more familiar $\normal_p\{0,J^{-1}\}$ result. 


\smallskip
{\csc Example.} 
In addition to the median $M_n$ in the situation of 3A, 
look at the mean absolute deviation statistic 
$\hatt\tau_n=n^{-1}\sumin |Y_i-M_n|$. 
We will show simultaneous convergence of 
$\sqrt{n}(M_n-\mu,\hatt\tau_n-\tau)$, where $\tau=\E|Y_i-\mu|$,
and for this assume finite variance of the $Y_i$'s. 

This can be accomplished by considering the parametric model
$f(y,\mu,\tau)=(2\tau)^{-1}\exp\{-|y-\mu|/\tau\}$ for data.
This model may be quite inadequate to describe the behaviour of
the data sequence, but the ML estimates are nevertheless
$M_n$ and $\hatt\tau_n$ as above. The traditional theorems
on ML behaviour require more smoothness than is present here,
and indeed often require that the true $f$ belongs to the model, 
but Theorem 4.1 can be used.
This is because $\log f(y,\mu,\tau)$ is concave in $(\mu,1/\tau)$. 
Verifying conditions involves details similar to those in 3A,
and we omit them here. The result is 
$$\pmatrix{\sqrt{n}(M_n-\mu) \cr
	   \sqrt{n}(\hatt\tau_n-\tau) \cr}
	   \arr_d\normal_2\bigl\{\pmatrix{0 \cr 0 \cr},
	   \pmatrix{1/\{4f(\mu)^2\}, &{\rm cov} \cr
	   	    {\rm cov}, & {\rm Var}\,Y_i-\tau^2 \cr}\bigr\}. $$ 
where the covariance is $\E I\{Y_i\le\mu\}|Y_i-\mu|-\half\tau$.
Note that there is asymptotic independence if $f$ is symmetric around $\mu$. 

\subsection
{\sl 4B. Bayes and maximum likelihood estimators 
are asymptotically equivalent.} 
It is well known that Bayes and ML estimation
are asymptotically equivalent procedures in regular situations.
In other words, if $\theta_n^*$ is the Bayes estimator under some
prior $\pi(\theta)$, then 
$\sqrt{n}(\theta_n^*-\theta_0)$ has the same 
limit distribution as $\sqrt{n}(\hatt\theta_n-\theta_0)$.
The standard proofs of this fact involve many technicalities, and
furthermore are typically restricted to calculations under the assumption that
the underlying $f(y,\theta_0)$  model is exactly correct, see e.g.~Lehmann
(1983, chapter 6.7).  
Below follows a reasonably quick proof of this fact,
and it is reassuring that the result is valid also
outside model circumstances. 


Let $\pi(\theta)$ be a prior density, assumed 
continuous at $\theta_0$ and satisfying the growth constraint
$$\pi(\theta) \le C_1 \exp(C_2|\theta|)\qt{for all }\theta, $$
where $C_1$ and $C_2$ are 
positive constants. The posterior density 
is proportional to $L_n(\theta)\pi(\theta)$,
where $L_n(\theta)=\prod_{i\le n}f(Y_i,\theta)$ is the likelihood. 
The Bayes estimator $\theta_n^*$ (under quadratic loss) 
is the posterior mean. Note that improper priors are accepted too. 

We shall make use of the following dominated convergence fact, 
which is a special case of Lemma A3 in the appendix.
Suppose $\{G_n(s,\omega)\}$ is a sequence of random
functions (assumed jointly measurable) such that
$G_n(s,\omega)\to G(s)$ in probability, for each $s$. 
Suppose $H(s)$ is an integrable function for which the set
$\{\omega\colon|G_n(s,\omega)|\le H(s)\hbox{ for all }s\}$
has probability tending to one. Then
$\int G_n(s,\omega)\,\d{}s\to\int G(s)\,\d{}s$ in probability. 
(Apply Lemma~A3 with~$X_n$ equal to~$G_n$ restricted to the set where $G_n\le
H$.)

\def\MLE{\hatt\theta_n}

{{\smallskip\sl
{\csc Theorem 4.2.} 
Under the conditions of Theorem 4.1, the  MLE estimator $\hatt\theta_n$ 
and the posterior mean~$\theta_n^*$ are 
asymptotically equivalent, in the sense that
$\sqrt{n}(\theta_n^*-\hatt\theta_n)\arr_p0$.
\smallskip}}

{\csc Proof:} 
Define the random convex function $A_n(s)$ by
$$\exp(-A_n(s)) = L_n(\hatt\theta_n+s/\sqrt{n})/L_n(\hatt\theta_n).$$
By definition of the ML estimator, 
$A_n$ achieves its minimum value of zero at~$s=0$.
By the change of variables
$\theta=\MLE+s/\sqrt{n}$  we find
$$\theta_n^*
={\int\theta L_n(\theta)\pi(\theta)\,\d{}\theta
   \over \int\ L_n(\theta)\pi(\theta)\,\d{}\theta} 
=\MLE+{1\over \sqrt{n}}
	{\int s\exp(-A_n(s))\pi(\MLE+s/\sqrt{n})\exp(-C_2|\MLE|)\,\d{}s
	\over \int \exp(-A_n(s))\pi(\MLE+s/\sqrt{n})\exp(-C_2|\MLE|)\,\d{}s}.$$
The random  function $A_n$ converges in probability uniformly on compact
sets to $\half s'Js$.  Define 
$\gamma_n=\inf_{|t|=1}A_n(t)$. It converges in probability to
$\gamma_0=\inf_{|t|=1}\half t'Jt>0$. 
Argue as in Lemma~2 to show that $A_n(s)\ge\gamma_n|s|$ for $|s|>1$. 
The domination condition needed for the fact noted above holds in both
numerator and denominator with
$$H(s)=\cases{2C_1 &if $|s|\le 1$, \cr
C_1|s|\exp(-\half\gamma_0|s|) &if $|s|>1$. \cr}$$
The ratio of integrals converges in probability to
$${\int s\exp(-\half s'Js)\,\pi(\theta_0)\exp(-C_2|\theta_0|)\,\d{}s
\over
\int \exp(-\half s'Js)\,\pi(\theta_0)\exp(-C_2|\theta_0|)\,\d{}s} = 0.$$
The result follows.
\square

\bigskip
{\bf 5. Logistic regression.} 
Suppose that $p+1$-tuplets $(x_i,Y_i)$ are observed,
where $x_i=(x_{i,1},\ldots,x_{i,p})'$ is a covariate vector 
`explaining' the binomial outcome $Y_i$. The logistic 
regression model postulates that the $Y_i$'s are independent with 
$${\rm Pr}\{Y_i=1|x_i\}=q(x_i,\beta)=
	{\exp(\beta'x_i)\over 1+\exp(\beta'x_i)} 
	\quad {\rm for\ some\ }\beta=\beta_0, \eqno(5.1)$$ 
and the ML estimator 
$\hatt\beta_n=(\hatt\beta_{n,1},\ldots,\hatt\beta_{n,p})'$ 
maximises the log-likelihood function 
$$\sumin \bigl[Y_i\log q(x_i,\beta)
	+(1-Y_i)\log\{1-q(x_i,\beta)\}\bigr] 
  =\sumin \bigl[Y_i\beta'x_i-\log\{1+\exp(\beta'x_i)\}\bigr]. $$
Of course the asymptotic normality of this estimator is 
well known and widely used, 
but precise sufficient conditions are not easy
to find in the literature. 

We will soon arrive at such,
employing results of 2B, which are applicable 
since the summands above are concave in $\beta$. 
As a preparatory exercise we mark down the following little expansion, 
which holds for all $u$ and $u+h$, 
in terms of $\pi(u)=\exp(u)/\{1+\exp(u)\}$: 
$$\log{1+\exp(u+h)\over 1+\exp(u)}
	=\pi(u)h+\half \pi(u)\{1-\pi(u)\}h^2
	 +\sixth\pi(u)\{1-\pi(u)\}\gamma(u,h)h^3, \eqno(5.2)$$
where $|\gamma(u,h)|\le \exp(|h|)$. This is proved from 
the exact third order Taylor expansion expression, with 
third term equal to $\sixth\pi(u')\{1-\pi(u')\}\{1-2\pi(u')\}h^3$, 
for appropriate $u'$ between $u$ and $u+h$. 
Some analysis reveals that 
$\pi(u')\{1-\pi(u')\}\le \exp(|h|)\,\pi(u)\{1-\pi(u)\}$,
regardless of $u$ and $h$. 
This is in fact quite similar 
to what results from using Lemma A2 in the appendix, but the bound 
on the remainder obtained here suits the problem better. 


\subsection
{\sl 5A. Under model conditions.}
In the spirit of our two aims, laid out in the Introduction,
we will first give a simpler result with a `pedagogical proof',
and then sharpen the tools to reach a second result with 
minimal regularity conditions. Under model conditions (5.1),
write for convenience $q_i=q(x_i,\beta_0)$, and let 
$J_n=\sumin q_i(1-q_i)x_ix_i'$ be the information matrix.  

{{\smallskip\sl
{\csc Theorem 5.1.} 
Assume that $\mu_n=\maxin|x_i|/\sqrt{n}\arr0$ 
and that $J_n/n\arr J$. Then, under model conditions (5.1),  
$\sqrt{n}(\hatt\beta_n-\beta_0)\arr_d\normal_p\{0,J^{-1}\}$.
\smallskip}}

{\csc Proof:} 
We will use Theorem 2.2 with 
$g_i(y_i,\beta)=\log f_i(y_i,\beta)=y_i\beta'x_i-\log\{1+\exp(\beta'x_i)\}$. 
The expansion noted above yields 
$$\eqalign{
\log {f_i(y_i,\beta_0+t)\over f_i(y_i,\beta_0)}
&=y_it'x_i-\bigl[\log\{1+\exp(\beta_0'x_i+t'x_i)\}
	  -\log\{1+\exp(\beta_0'x_i)\}\bigr] \cr
&=(y_i-q_i)x_i't
	-\half q_i(1-q_i)(t'x_i)^2
	-\sixth q_i(1-q_i)\gamma_i(t)(t'x_i)^3 \cr
&=D_i(y_i)'t-R_i(y_i,t). \cr}$$
Here $D_i(y_i)=(y_i-q_i)x_i$ and 
$R_i(y_i,t)=\half t'q_i(1-q_i)x_ix_i't+v_{i,0}(t)$,
where $|\gamma_i(t)|\le\exp(|t'x_i|)$ in the expression for $v_{i,0}(t)$.
Note that $J_n=K_n$, in the notation of Theorem 2.2, 
and that $R_i(Y_i,t)$ has 
zero variance, so what we have to prove is 
(i) that $\sumin v_{i,0}(s/\sqrt{n})\arr0$,
(ii) that the Lindeberg conditions are satisfied for 
$\sumin n^{-1/2}(Y_i-q_i)x_i$.
But 
$$\eqalign{
\Big|\sumin v_{i,0}(s/\sqrt{n})\Bigr|
	&\le\sumin \sixth q_i(1-q_i)
	\exp(|s'x_i/\sqrt{n}|)\,|s'x_i/\sqrt{n}|^3 \cr
	&\le\sumin \sixth q_i(1-q_i)\exp(|s|\mu_n)
		\,(s'x_ix_i's/n)\,|s|\mu_n \cr
	&=\sixth|s|\mu_n\exp(|s|\mu_n)\,s'(J_n/n)s, \cr}$$
which goes to zero. 
And the Lindeberg condition is that for each $s$ and $\delta$
$$\sumin En^{-1}(Y_i-q_i)^2(s'x_i)^2\,I\{|(Y_i-q_i)s'x_i/\sqrt{n}|\ge\delta\}
	\arr0, $$
and this sum is bounded by $s'(J_n/n)s\,I\{|s|\mu_n\ge\delta\}$.
This ends the proof. \square

\smallskip
If the $x_i$'s are i.i.d.~from some covariate distribution $H$,
then $\mu_n\arr0$ a.s.~exactly when the components of $x_i$ have 
finite second moment. 
This also secures convergence of 
$J_n/n$ to $J=\int q(x,\beta_0)\{1-q(x,\beta_0)\}\,xx'H(\d x)$. 

Our second and sharper theorem is proved next, by 
squeezing more out of the bound of the $v_{i,0}(t)$ 
remainder and more out of the Lindeberg condition.

{{\smallskip\sl 
{\csc Theorem 5.2.} 
Assume that the $\lambda_n=\maxin|J_n^{-1/2}x_i|$
sequence is bounded, and that  
$$N_n(\delta)=\sumin q_i(1-q_i)\,x_i'J_n^{-1}x_i\,
	I\{|J_n^{-1/2}x_i|\ge\delta\}
	\arr0 \quad {\sl for\ each\ positive\ }\delta. \eqno(5.3)$$
Then, under model conditions (5.1), 
$J_n^{1/2}(\hatt\beta_n-\beta_0)\arr_d\normal_p\{0,I_p\}$.	 		
\smallskip}}

{\csc Proof:}
We consider the random convex function 
$\sumin\{\log f_i(Y_i,\beta_0+J_n^{-1/2}s)-\log f_i(Y_i,\beta_0)\}$,
which upon using the expansion again can be rearranged as
$U_n's-\half s's-r_n(s)$, where 
$U_n=J_n^{-1/2}\sumin(Y_i-q_i)x_i$ and 
$r_n(s)=\sumin\sixth q_i(1-q_i)
	\gamma_i(s'J_n^{-1/2}x_i)(s'J_n^{-1/2}x_i)^3$.
We are to prove (i) that $r_n(s)\arr0$, 
and (ii) that $U_n\arr_d\normal_p\{0,I_p\}$. 

At this stage we call on appendix A1
where it is shown that (5.3) is a sufficient 
and actually also a necessary condition (ii) to hold. 
And $|r_n(s)|$ is bounded by 
$\sumin\sixth q_i(1-q_i)\exp(|s'J_n^{-1/2}x_i|)|s'\allowbreak J_n^{-1/2}x_i|^3$.
We split this sum into $|J_n^{-1/2}x_i|<\delta$ 
summands and $|J_n^{-1/2}x_i|\ge\delta$ summands. 
The first sum is bounded by 
$\sixth|s|^3\delta \exp(|s|\delta)$,
and the second is bounded by 
$\sixth|s|^3\lambda_n\exp(|s|\lambda_n)\,N_n(\delta)$. 
Letting $n\arr\infty$ and $\delta\arr0$ afterwards
shows that indeed $r_n(s)\arr0$. \square 


\smallskip
It is worth noting that the $N_n(\delta)\arr0$ condition in the theorem 
serves two purposes: forcing an analytic remainder term towards zero,
and securing uniform neglibility of individual terms in the large-sample 
distribution of $J_n^{-1/2}\sumin D_i(Y_i)$, i.e.~a normal limit.
Note also that $\lambda_n\arr0$ suffices for the conclusion to
hold, since $N_n(\delta)\le p\lambda_n/\delta$. 

\subsection
{\sl 5B. Outside model conditions.}
Let us next depart from the strict model assumption (5.1),
which in most cases merely is intended to provide a 
reasonable approximation to some more complicated reality,
and stipulate only that ${\rm Pr}\{Y=1|x\}=q(x)$ for some
true, underlying $q(x)$ function. Fitting the logistic regression
equation makes sense still, and turns out to aim 
at achieving the best approximation $q(x,\beta)$ to the true $q(x)$, 
in a sense made precise as follows. Let 
$$\Delta_x[q(x),q(x,\beta)]=q(x)\log {q(x)\over q(x,\beta)}
	+\{1-q(x)\}\log{1-q(x)\over 1-q(x,\beta)} $$
be the Kullback--Leibler distance from true binomial 
$(1,q(x))$ to modelled binomial $(1,q(x,\beta))$, and let 
$\Delta[q(.),q(.,\beta)]=\int \Delta_x[q(x),q(x,\beta)]\,H(\d x)$
be the weighted distance between the true probability curve to
the modelled probability curve, in which $H$ again is the  
`covariate distribution' for $x$'s, as discussed in 2C. 
The following can now be proved using methods of 2C: 
ML estimation is $\sqrt{n}$-consistent for the value $\beta_0$
that minimises the weighted Kullback--Leibler distance $\Delta$,
and $\sqrt{n}(\hatt\beta_n-\beta_0)\arr_d\normal_d\{0,J^{-1}KJ^{-1}\}$,
provided the two matrices  
$$J=\E\,XX'q(X,\beta_0)\{1-q(X,\beta_0)\}
	=\int xx'\,q(x,\beta_0)\{1-q(x,\beta_0)\}\,H(\d x), $$
$$K=\E\,XX'\{Y-q(X,\beta_0)\}^2
	=\int xx'\bigl[q(x)\{1-q(x)\}+\{q(x)-q(x,\beta_0)\}^2\bigr]\,H(\d x) $$
are finite. This result was also obtained in Hjort (1988a),
where various implications for statistical inference also are discussed. 

\bigskip
{\bf 6. Cox regression.}
In this section new proofs are 
presented for the consistency and asymptotic normality 
of the usual estimators in Cox's famous semiparametric  
regression model for survival analysis data. 
The parametric Cox regression model is somewhat simpler,
and is treated in 7A below.  
The regularity requirements we need turn out in both
cases to be weaker than those earlier presented in the literature. 


The most complete results and proofs in the literature 
for the basic large-sample properties 
of the estimators in this model are perhaps those of 
Andersen and Gill (1982) and Hjort (1992). 
Andersen and Gill obtain results under the conditions of
the model, and with regularity conditions quite weaker than 
earlier i.i.d.~type assumptions, whereas Hjort explores 
the large-sample behaviour also outside the conditions of the model. 
For a history of the Cox model 
and the various approaches to reach asymptotics results, see 
Andersen, Borgan, Gill \& Keiding (1992, chapter VII).

Our present intention is to provide yet another proof, 
which in several ways is simpler and requires less
involvement with the martingale techniques 
than the one of Andersen and Gill.
As in the previous section we choose to present 
two theorems, reflecting our two aims explained in section 1. 
The first holds when the covariates are bounded,
in which case the proof is quite transparent,
and extra regularity conditions can be kept quite minimal.
The second version is more sophisticated in that it tolerates
unbounded covariates and weakens regularity conditions further. 

The usual Cox regression model for possibly 
censored lifetimes with covariate information 
is as follows: The individuals have independent lifetimes
$T_1^0,\ldots,T_n^0$, and the $i$-th has hazard rate
$$\lambda_i(s)=\lambda(s)\exp(\beta'z_i(s))
  =\lambda(s)\exp(\beta_1z_{i,1}(s)+\cdots \beta_pz_{i,p}(s)), \eqno(6.1)$$
depending on that person's covariate vector 
$z_i(s)$, and involving some unspecified basis hazard rate $\lambda(s)$. 
As indicated the covariates are allowed to depend on time $s$,
and they can be random processes, as long as they are previsible; 
$z_i(s)$ should only depend on information available at time $s-$
(for a full discussion of previsibility, or predictability,
see Andersen et al.~(1992, p.~65--66)). 
There is a possibly interfering censoring time $C_i$ leaving just 
$T_i=\min\{T_i^0,\allowbreak C_i\}$ and $\delta_i=I\{T_i^0\le C_i\}$
to the statistician. 
Consider the at risk indicator function $Y_i(s)=I\{T_i\ge s\}$,
which is left continuous and hence previsible, 
and the counting process $N_i$ with
mass $\delta_i$ at $T_i$, i.e.~$\d N_i(s)=I\{T_i\in[s,s+ds],\delta_i=1\}$.
The log partial likelihood can then be written 
$$G_n(\beta)=\sumin \int_0^L\bigl\{\beta'z_i(s)
	-\log R_n(s,\beta)\bigr\}\,\d{}N_i(s), \eqno(6.2) $$ 
featuring the empirical relative risk function 
$R_n(s,\beta)=\sumin Y_i(s)\exp(\beta'z_i(s))$;
see for example Andersen et al.~(1992, chapter VII).
It is assumed that data are collected on the finite
time interval $[0,L]$ only. 
The {\it Cox estimator} is the value $\hatt\beta_n$ 
that maximises the partial likelihood.


Lemma A2 of the appendix allows us an expansion for 
$\log R_n(s,\beta_0+x)$,
using $w_i=Y_i(s)\exp(\beta_0'\allowbreak z_i(s))$ and $a_i=z_i(s)'x$.
The result is 
$$\log R_n(s,\beta_0+x)-\log R_n(s,\beta_0)
	=\bar z_n(s)'x+\half x'V_n(s)x+v_n(x,s), \eqno(6.3)$$
where 
$$\bar z_n(s)=\sumin p_{n,i}(s)z_i(s)
	\quad {\rm and} \quad 
V_n(s)=\sumin p_{n,i}(s)(z_i(s)-\bar z_n(s))
	(z_i(s)-\bar z_n(s))', \eqno(6.4)$$
and $p_{n,i}(s)=Y_i(s)\exp(\beta_0'z_i(s))/R_n(s,\beta_0)$.
A bound for the remainder term in (6.3) is 
$|v_n(x,s)|\le\hbox{$4\over3$}\maxin|(z_i(s)-\bar z_n(s))'x|^3.$
Observe that $\bar z_n(s)$ and $V_n(s)$ can be interpreted
as the mean value and the variance matrix for $z_i(s)$,
where this covariate vector is randomly selected among 
those at risk at time $s$ with probabilities proportional to 
the relative risks $\exp(\beta_0'z_i(s))$. 

All this leaves us suitably prepared for a theorem. 

{\smallskip
{\csc Theorem 6.1.} \sl
Assume that the hazard rate for the $i$'th individual follows 
the Cox model (6.1) with a true parameter $\beta_0$ and 
a continuous positive basis hazard $\lambda(s)$, 
and that the covariate processes $z_i(s)$ are previsible and
uniformly bounded.   
Assume that 
$$J_n(s)=n^{-1}\sumin Y_i(s)\exp(\beta_0'z_i(s))\,
	(z_i(s)-\bar z_n(s))(z_i(s)-\bar z_n(s))'\arr_pJ(s)
	\eqno(6.5)$$
for almost all $s$ in $[0,L]$ 
and that $J=\int_0^LJ(s)\lambda(s)\,\d{}s$ is positive definite.
Then $\hatt\beta_n$ is $\sqrt{n}$-consistent for $\beta_0$ and
$\sqrt{n}(\hatt\beta_n-\beta_0)\arr_d \normal_p\{0,J^{-1}\}$. 
\smallskip}

{\csc Proof:}
As a simple consequence of earlier efforts we have 
$$\eqalign{
G_n^*(x)&=G_n(\beta_0+x/\sqrt{n})-G_n(\beta_0) \cr
	&=\sumin \int_0^L \bigl[n^{-1/2}(z_i(s)-\bar z_n(s))'x
	-\half n^{-1}x'V_n(s)x-v_n(x/\sqrt{n},s)\bigr]\,\d{}N_i(s) \cr
	&=U_n'x-\half x'J_n^* x-r_n(x), \cr} \eqno(6.6)$$	
where we write 
$$U_n=n^{-1/2}\sumin \int_0^L(z_i(s)-\bar z_n(s))\,\d{}N_i(s) 
	\quad {\rm and} \quad 
  J_n^*=n^{-1}\int_0^L V_n(s)\,\d{}\bar N_n(s), \eqno(6.7)$$
using $\bar N_n(.)=\sumin N_i(.)$ to denote 
the aggregated counting process for the data. 
The designated remainder term is 
$r_n(x)=\int_0^Lv_n(x/\sqrt{n},s)\,\d{}\bar N_n(s)$, 
which goes to zero, since it is bounded by 
$\int_0^L{4\over3}(2K)^3|x|^3/n^{3/2}\,\d{}\bar N_n(s)$,
which is $O(n^{-1/2})$. The $K$ here is the absolute bound 
on the covariates. 
That the (6.6) function is concave in $x$ 
is clear from the convexity of $\log R_n(s,\beta)$ in $\beta$. 
By the basic method of section 1 it only remains to show 
(i) that $J_n^*\arr_pJ$ and 
(ii) that $U_n\arr_d \normal_p\{0,J\}$. 

At this stage we need some of the easier bits of the
martingale representation and convergence theory 
for counting processes, but manage to avoid needing 
some of the more sophisticated inequalities and technicalities 
that have invariably been present in earlier rigorous
proofs, like in Andersen and Gill (1982).
The counting process $N_i$ has compensator process 
$A_i(t)=\int_0^t Y_i(s)\exp(\beta_0'z_i(s))\,\d{}\Lambda(s)$,
writing $\d\Lambda(s)=\lambda(s)\,\d{}s$. This means that 
$M_i(t)=N_i(t)-A_i(t)$ is a zero mean martingale,
with increments 
$\d M_i(s)=\d N_i(s)-Y_i(s)\exp(\beta_0'z_i(s))\,\d{}\Lambda(s)$.   
One can show that $M_i(t)^2-A_i(t)$ 
as well as $M_i(t)M_j(t)$ are martingales too, when $i\not=j$, 
which in martingale theory parlance 
means that $M_i$ has variance process 
$\langle M_i,M_i\rangle(t)=A_i(t)$
and that they are orthogonal, i.e.~$\langle M_i,M_j\rangle=0$. 
See Andersen et al.~(1992, chapter II), for example. 
Inserting $\d N_i(s)=\d M_i(s)+dA_i(s)$ in (6.7) 
leads to 
$$J_n^*=\int_0^LJ_n(s)\,\d{}\Lambda(s)
	+n^{-1}\sumin \int_0^LV_n(s)\,\d{}M_i(s) \eqno(6.8)$$
and 
$$U_n=n^{-1/2}\sumin \int_0^L(z_i(s)-\bar z_n(s))\,\d{}M_i(s), 
	\eqno(6.9)$$
in that two other terms cancel. 

We are now in a position to prove (i) and (ii). 
Note that the first term of (6.8) goes to $J$ in probability
by boundedness of the integrand and Lemma A3 in the appendix. 
The second term is $O_p(n^{-1/2})$,
which can be seen using boundedness of covariates 
in conjunction with the result 
$$\E\Bigl\{\int_0^L\sumin H_i(s)\,\d{}M_i(s)\Bigr\}^2
	=\E\sumin\int_0^LH_i(s)^2\,\d{}\langle M_i,M_i\rangle(s), $$
valid for previsible random functions $H_i$. 
This proves (i). 
To prove convergence in distribution of $U_n$ we essentially use 
the version of Rebolledo's martingale central limit theorem 
given in Andersen and Gill (1982, appendix I). 
Its variance process converges properly, 
$$\langle U_n,U_n\rangle(L)
	=n^{-1}\sumin \int_0^L
	(z_i(s)-\bar z_n(s))(z_i(s)-\bar z_n(s))'
		\,\d{}\langle M_i,M_i\rangle(s) 
	=\int_0^L J_n(s)\,\d{}\Lambda(s)\arr_p J, $$
and the necessary Lindeberg type condition is also satisfied:
$$n^{-1}\sumin \int_0^L|z_i(s)-\bar z_n(s)|^2\,
	I\{n^{-1/2}|z_i(s)-\bar z_n(s)|\ge\eps\}\,
	Y_i(s)\exp(\beta_0'z_i(s))\,\d{}\Lambda(s) \arr_p0 \eqno(6.10)$$
since the indicator function ends up being zero for all
large $n$. \square
	
\smallskip
Next we present a stronger theorem with weaker conditions imposed.
The proof is basically the same as for the previous result,
but more is squeezed out of bounds for remainder terms and 
out of conditions for the martingale convergence to hold. 

{\smallskip\sl
{\csc Theorem 6.2.} 
Assume that the hazard rate for the $i$'th individual follows 
the Cox model (6.1) with a true parameter $\beta_0$ and 
a continuous positive basis hazard $\lambda(s)$.
Assume that $J_n(s)$ goes to some $J(s)$ 
in probability for almost all $s$, as in (6.5),  
and that $\int_0^LJ_n(s)\lambda(s)\,\d{}s\arr_pJ=\int_0^LJ(s)\lambda(s)\,\d{}s$,
a positive definite matrix. Suppose finally that 
$$\mu_n(s)=n^{-1/2}\maxin|z_i(s)-\bar z_n(s)|\arr_p0 
	\quad {\sl for\ almost\ each\ }s \eqno(6.11)$$
and that $\max_{s\le L}\mu_n(s)$ is stochastically bounded. Then again 
$\sqrt{n}(\hatt\beta_n-\beta_0)\arr_d \normal_p\{0,J^{-1}\}$. 
\smallskip}

{\csc Proof:}
(6.6) and (6.7) still hold, and we plan to demonstrate 
(i) $r_n(x)\arr_p0$,
(ii) $J_n^*\arr_pJ$, and
(iii) $U_n\arr_d\normal_p\{0,J\}$.

(i) is proved by using the tighter bound for $v_n(x,s)$ of (6.3) 
available by employing Lemma A2, namely 
${2\over3}g(\maxin|(z_i(s)-\bar z_n(s))'x|)\,x'V_n(s)x$,
for $g(u)=u\exp(2u+4u^2)$. This leads to 
$$|r_n(x)|\le\int_0^L\hbox{$2\over3$}g(\mu_n(s)|x|)\,x'V_n(s)x 
	\,\d{}\bar N_n(s)/n.$$ 
Split this into two terms, using 
$\d \bar N_n(s)=R_n(s,\beta_0)\,\d{}\Lambda(s)+\sumin \d M_i(s)$.
The first of the resulting terms goes to zero in probability 
by assumptions on $J_n(s)$ and dominated convergence (appendix A3),
and the other term is of smaller stochastic order.
Secondly (ii) follows as in the previous proof, since 
the second term of (6.8) vanishes in probability, 
by variations of the same arguments.
Finally two ingredients are needed to secure (iii). 
The first is $\langle U_n,U_n\rangle(L)\arr_pJ$, which 
holds by assumptions as in the previous proof. 
The second is a more elaborate demonstration 
of the Lindeberg type condition (6.10), now accomplished by bounding
it with		
$$\int_0^L{\rm Tr}(J_n(s))\,I\{\mu_n(s)\ge\eps\}\,\d{}\Lambda(s), $$
which goes to zero in probability 
by dominated convergence (the integrand goes pointwise to zero
in probability and is dominated by ${\rm Tr}(J_n(s))$, 
see appendix A3 again). 

And all this combined with the Basic Corollary triumphantly implies 
that the argmax of the (6.6) function, which is 
$\sqrt{n}(\hatt\beta_n-\beta_0)$, is only $o_p(1)$ away from 
the argmax of $U_n'x-\half x'Jx$, which is $J^{-1}U_n$.
This proves consistency and asymptotic normality. \square  

\smallskip
{\csc Remarks.}
(i) 
Usually one would have $V_n(s)\arr_pV(s)$ and 
$n^{-1}R_n(s,\beta_0)\arr_pR(s,\beta_0)$, say,
so that $J_n(s)=V_n(s)R_n(s,\beta_0)/n\arr_p J(s)=V(s)R(s,\beta_0)$;
in particular $J=\int_0^LV(s)R(s,\beta_0)\,\d{}\Lambda(s)$ in such cases,
and this is the expression typically encountered for the inverse
covariance matrix. 
(ii) 
The Andersen and Gill regularity requirements include rather strong 
uniform convergence statements, in both time $s$ and $\beta$ near $\beta_0$.
In the development above this would mean requiring 
$$\sup_{s\in[0,L]}\sup_{\beta\in U(\beta_0)}
\Big|n^{-1}\sumin Y_i(s)z_i(s)z_i(s)'\exp(\beta'z_i(s))
	-J(s,\beta)\Big|\arr_p0, $$
for example, for a suitable neighbourhood $U(\beta_0)$ and 
a suitable limit function $J(s,\beta)$. This contrasts sharply
with our condition (6.5), which is only about $\beta_0$,
and is pointwise in $s$. Andersen and Gill also include 
various other asymptotic stability conditions, about uniform
continuity and differentiability in $\beta$ of their limit
functions, that are not needed here. Similarly, their conditions
almost require $\max_{s\le L}\mu_n(s)\arr_p0$ where we come
away with pointwise convergence.
(iii) 
is interesting to see that the key requirement 
(6.11) serves two different purposes: forcing an analytical remainder
term towards zero as well as securing uniform negligibility 
of individual terms, i.e.~limiting normality. 
(iv) 
The methods used here can be applied to solve the large-sample
behaviour problem also outside model conditions, say when
the true hazard rate is $\lambda(s)\,r(z_{i,1}(s),\ldots,z_{i,p}(s))$  
for individual $i$. See Hjort (1992) for results. 
There are also various alternative estimation techniques that
can be employed in the Cox model, see for examples Hjort (1991)
for local likelihood smoothing and Hjort (1992) for weighted log partial
likelihood estimation. Again techniques from the present paper 
can be applied. 
(v) 
Finally Jeffreys type arguments can be given in favour of using 
the vague prior $\pi(\beta)=1$, see Hjort (1986),
where it is also shown that the (improper) pseudo-Bayes estimator 
$\beta_n^*=\int\beta\exp(G_n(\beta))\,\d \beta/\int\exp(G_n(\beta))\,\d \beta$
is asymptotically equivalent to the Cox estimator $\hatt\beta_n$. 
The arguments of 4B can be used to provide a quicker and simpler proof of this.

\bigskip
{\bf 7. Further applications.}

\subsection
{\sl 7A. Exponential hazard rate regression.} 
The traditional Cox model (6.1) is semiparametric, since 
the basis hazard rate $\lambda(.)$ there is left unspecified.
The parametric version 
$$\lambda_i(s)=\lambda_0(s)\,\exp(\beta'z_i(s))
	=\lambda_0(s)\,\exp(\beta_1z_{i,1}(s)+\cdots+\beta_pz_{i,p}(s)),
		\eqno(7.1) $$
where $\lambda_0(.)$ is fully specified (equal to 1, for example), 
is also important in survival data analysis. 
Let us briefly show that the arguments above efficiently 
lead to a precise theorem about the maximum likelihood estimator 
in this model as well. 

With notation and assumptions otherwise being 
as in section 6 the log-likelihood can be written 
$$\log L_n(\beta)=\sumin \int_0^L\bigl\{\beta'z_i(s)\,\d{}N_i(s)
	-Y_i(s)\exp(\beta'z_i(s))\,\d{}s\bigr\}, $$
see for example Andersen et al.~(1992, chapter VI), 
and let $\hatt\beta_n$ be 
the ML estimator maximising this expression.  	
Assume that data follow (7.1) for a certain $\beta_0$. 
Using martingales 
$\d M_i(s)=\d N_i(s)-Y_i(s)\exp(\beta_0'\allowbreak z_i(s))\,\d{}\Lambda_0(s)$,
writing $\d\Lambda_0(s)=\lambda_0(s)\,\d{}s$, we find 
$$\eqalign{G_n(x)&=\log L_n(\beta_0+x)-\log L_n(\beta_0) \cr
&=\sumin	\int_0^L\Bigl[x'z_i(s)\,\bigl\{\d M_i(s)
	+Y_i(s)\,\exp(\beta_0'z_i(s))\,\d{}\Lambda_0(s)\bigr\} \cr
&\qquad\qquad\qquad 
 	-Y_i(s)\,\exp(\beta_0'z_i(s))\,\{\exp(x'z_i(s))-1\}\,\d{}s\Bigr] \cr
&=\Bigl(\sumin\int_0^L\!\! z_i(s)\,\d{}M_i(s)\Bigr)'x 
	-\sumin\int_0^L\!\! Y_i(s)\,\exp(\beta_0'z_i(s))\,
	\bigl\{\exp(x'z_i(s))-1-x'z_i(s)\bigr\}\,\d{}\Lambda_0(s). \cr}$$
Now use $|\exp(u)-1-u-\half u^2|\le\sixth|u|^3\exp(|u|)$, and introduce 
$$J_n=\sumin\int_0^LY_i(s)\exp(\beta_0'z_i(s))\,z_i(s)z_i(s)'\,\d{}\Lambda_0(s) 
	{\rm \ and\ } 
  U_n=J_n^{-1/2}\sumin\int_0^Lz_i(s)\,\d{}M_i(s) \eqno(7.2)$$
to find 
$G_n(J_n^{-1/2}x)=U_n'x-\half|x|^2-r_n(x)$, 
with a remainder bound 
$$|r_n(x)|\le \sixth\sumin\int_0^L
	Y_i(s)\exp(\beta_0'z_i(s)+|x'J_n^{-1/2}z_i(s)|)
	\sixth|x'J_n^{-1/2}z_i(s)|^3\,\d{}s. \eqno(7.3)$$
To formulate a theorem with quite weak conditions, let
$J_n(s)=\sumin Y_i(s)\exp(\beta_0'z_i(s))\,z_i(s)z_i(s)'$, 
so that the `observed information matrix' is 
$J_n=\int_0^LJ_n(s)\,\d{}\Lambda_0(s)$. 

{{\smallskip\sl
{\csc Theorem 7.1.}
Let $\lambda_0(s)$ be positive and continuous on $[0,L]$. 
Suppose there is a $c_n$ sequence converging to infinity such that 
$J_n(s)/c_n$ for almost all $s$ goes in probability to some $J(s)$,
and that $J_n/c_n\arr_p J=\int_0^LJ(s)\,\d{}\Lambda_0(s)$,
where this limit matrix is positive definite. 
Assume furthermore that for almost all $s$, 
$$N_n(s,\delta)=\sumin 
	z_i(s)'J_n^{-1}z_i(s)\,Y_i(s)\exp(\beta_0'z_i(s))\,
 	I\{|J_n^{-1/2}z_i(s)|\ge\delta\}\arr_p0 \eqno(7.4)$$
for each $\delta>0$, and that $\mu_n(s)=\maxin|J_n^{-1/2}z_i(s)|$ 
is stochastically bounded, uniformly in $s$. 
Then $J_n^{1/2}(\hatt\beta_n-\beta_0)\arr_d\normal_p\{0,I_p\}$. 
\smallskip}}

Some brief remarks are in order before turning to the proof. 
(i) $z_i(s)J_n^{-1}z_i(s)$ can be replaced by $z_i(s)'J^{-1}z_i(s)/c_n$
here, and $|J_n^{-1/2}z_i(s)|$ with $|J^{-1/2}z_i(s)|/c_n^{1/2}$.
(ii) In many practical situations the $c_n$ will be equal to $n$. 
(iii) The elements of $J_n$ may in some cases concievably go to
infinity with different rates, and then the `asymptotic stability'
requirement should be the existence of matrices $C_n$ going 
to infinity such that $C_n^{-1}J_n(s)\arr J(s)$ etcetera. 
The theorem still holds. 
(iv) In many cases one would have $\mu_n(s)\arr_p 0$ for almost all $s$,
and this implies condition (7.4),
since in fact $D_n(s,\delta)\le pI\{\mu_n(s)\ge\delta\}$. 
(v) If the $z_i(s)$ covariate processes are uniformly 
bounded, then (iv) applies and hence the conclusion.
(vi) Our conditions are much weaker than those used 
elsewhere to secure large sample normality, see 
for example Borgan (1984, section 6). 
(vii) Finally we note that the proof below becomes easier under 
circumstances (iv) or (v).  

\smallskip
{\csc Proof:}
The log-likelihood is concave by Lemma A.2 and hence so is 
the $G_n(J_n^{-1/2}x)$ function. We are to prove 
(i) $r_n(x)\arr 0$ in probability for each $x$, and  
(ii) that $U_n\arr \normal_p\{0,I_p\}$ in distribution. 

To prove (i) let $r_n(x,s)$ be the integrand in the bound 
occurring in (7.3), 
so that $|r_n(x)|\le\sixth\int_0^Lr_n(x,s)\,\d{}\Lambda_0(s)$. 
It will suffice to show that $r_n(x,s)\arr 0$ in probability 
for almost all $s$ and to bound it properly. Splitting into 
$|J_n^{-1/2}z_i(s)|<\delta$ terms and 
$|J_n^{-1/2}z_i(s)|\ge\delta$ terms we find 
$r_n(x,s)\le|x|^3\delta\exp(|x|\delta)
	+|x|^3\mu_n(s)\exp(|x|\mu_n(s))\,N_n(s,\delta)$,
after which the claim follows by our precautions and by 
the dominated convergence lemma of the appendix. 

Next (ii) can be replaced by $U_n^*\arr_d \normal_p\{0,I_p\}$, 
where 
$U_n^*=\sumin\int_0^L c_n^{-1/2}J^{-1/2}z_i(s)\,\d{}M_i(s)$, 
and we show this employing the Rebolledo theorem version given in 
Andersen and Gill (1982, appendix I). 
Its variance process converges properly, 
$$\eqalign{
\langle U_n^*,U_n^*\rangle(L)
	&=\sumin\int_0^Lc_n^{-1}J^{-1/2}z_i(s)z_i(s)'J^{-1/2}\,
		Y_i(s)\exp(\beta_0'z_i(s))\,\d{}\Lambda_0(s) \cr
	&=J^{-1/2}(J_n/c_n)J^{-1/2}\arr_pI_p, \cr}$$
and the Lindebergian condition is also satisfied:
$$\sumin\int_0^Lc_n^{-1}|J^{-1/2}z_i(s)|^2\,
	I\{c_n^{-1/2}|J_n^{-1/2}z_i(s)|\ge\delta\}\,
	Y_i(s)\exp(\beta_0'z_i(s))\,\d{}\Lambda_0(s)\arr_p0.$$
This holds since the integrand is asymptotically the same as 
$N_n(s,\delta)$ of (7.4), and is bounded by the constant $p$,
so that the dominated convergence lemma applies. \square

\subsection
{\sl 7B. Poisson regression.} 
Suppose $Y_1,\ldots,Y_n$ are independent counts with 
$$Y_i\sim{\rm Poisson}({\rm mean}_i), 
	\quad {\rm with\ }{\rm mean}_i=\exp(\beta'z_i), \eqno(7.5)$$
depending on a certain $p$-dimensional
covariate vector $z_i$, for a certain true parameter value $\beta_0$. 
It is convenient now to write $\exp(u)=1+u+\half u^2+\sixth\rho(u)$,
where a bound for the remainder function is $|\rho(u)|\le|u|^3\exp(|u|)$.
The log-likelihood 
$\log L_n(\beta)=\sumin \{Y_i\beta'z_i-\exp(\beta'z_i)\}$ 
is concave, and after development very similar to that of section 7A
one finds 
$$\log L_n(\beta_0+x)-\log L_n(\beta_0)
	=\Bigl(\sumin(Y_i-\mu_i)z_i\Bigr)'x
	-\half x'J_nx-\sixth v_n(x), $$
in which $J_n=\sumin \mu_iz_iz_i'$ and 
$v_n(x)=\sumin\mu_i\rho(x'z_i)$. In these expressions 
$\mu_i=\exp(\beta_0'z_i)$ is the true mean for $Y_i$ under the model. 

Before passing to a theorem we solve a relevant exercise 
in asymptotics of linear combinations of independent Poisson variables. 
If $Y_{n,i}$ is Poisson with mean $\mu_{n,i}$, 
then $\sumin(Y_{n,i}-\mu_{n,i})x_{n,i}$, normed such that its variance
$\sumin\mu_{n,i}x_{n,i}^2=1$, goes to a standard normal if and only if
$\sumin\mu_{n,i}\rho(tx_{n,i})\arr0$ for each $t$, which is equivalent to
$\sumin\mu_{n,i}\rho(|x_{n,i}|)\arr0$. This is seen after considering 
moment or cumulant generating functions. 

{{\smallskip\sl
{\csc Theorem 7.2.} 
Let $\hatt\beta_n$ be the ML estimator based on the first 
$n$ Poisson counts, assumed to follow (7.5) for 
a certain $\beta_0$, with means $\mu_i=\exp(\beta_0'z_i)$. Then 
$J_n^{1/2}(\hatt\beta_n-\beta_0)\arr_d \normal_p\{0,I_p\}$
if and only if $\sumin\mu_i\rho(|J_n^{-1/2}z_i|)\arr0$. 
A simple sufficient condition for this to hold is that 
$\lambda_n=\maxin|J_n^{-1/2}z_i|$ is bounded and that 
$\sumin\mu_i|J_n^{-1/2}z_i|^3\arr0$; or, equivalently,
that $\lambda_n$ is bounded and that 
$$N_n(\delta)=\sumin\mu_iz_i'J_n^{-1}z_i\,I\{|J_n^{-1/2}z_i|\ge\delta\}
	\arr0\quad{\sl for\ each\ }\delta. $$
\smallskip}}

{\csc Proof:} 
The function 
$\log L_n(\beta_0+J_n^{-1/2}x)-\log L_n(\beta_0)$
is concave in $x$ and can be written 
$U_n'x-\half |x|^2-r_n(x)$, where 
$U_n=J_n^{-1/2}\sumin(Y_i-\mu_i)z_i$ and where 
$r_n(x)=\sumin\mu_i\rho(x'J_n^{-1/2}z_i)$. 
This is quite similar to the situation in 7A, 
and the maximiser $J_n^{1/2}(\hatt\beta_n-\beta_0)$ goes to 
a standard $p$-dimensional normal if and only if 
(i) $r_n(x)\arr0$ and (ii) $U_n\arr_d\normal_p\{0,I_p\}$.  
But using the result above in tandem with the Cram\'er--Wold theorem
one sees that $\sumin\mu_i\rho(|J_n^{-1/2}z_i|)\arr0$ 
is necessary and sufficient for (ii), and indeed also 
necessary and sufficient for (i). 
The other statements of the theorem follow 
>from $|\rho(u)|\le|u|^3\exp(|u|)$. \square

\smallskip
We note that $\lambda_n\arr0$ is clearly sufficient for the result to hold.

\subsection
{\sl 7C. Generalised linear models.}
Consider a situation with independent $Y_i$'s from densities of the form 
$f(y_i|\theta_i)=\exp\{(y_i\theta_i-b(\theta_i))/a(\phi)+c(y_i,\phi)\}$, 
and where $\theta_i$ is parametrised as a linear $x_i'\beta$. 
This is a generalised linear model with canonical link, 
see McCullagh and Nelder (1989). 
The likelihood in $\beta$ is log-concave,
and theorems about the large-sample behaviour of the ML estimator, 
under very weak regularity conditions, can be written down and proved
by the methods exemplified in sections 5 and 7A. 

\subsection
{\sl 7D. Pseudo-likelihood estimation in Markov chains.}
\fermat{too much?}Suppose 
$X_0,X_1,\ldots$ forms a Markov chain on the state
space $\{1,\ldots,k\}$. Instead of focusing on 
transition probabilities,  
consider direct modelling of $X_i$ given its neighbours,
say $X_{\dell i}=x_{\dell i}$. 
A very flexible and convenient class of models is described by
$$f_\beta(x_i|x_{\dell i})
	={\rm const.}\,\exp\{\alpha_i(x_i)+\beta'H_i(x_i,x_{\dell i})\}
	 ={\exp\{\alpha_i(x_i)+\beta'H_i(x_i,x_{\dell i})\}
	   \over
	   \sum_{j=1}^k\exp\{\alpha_i(j)+\beta'H_i(j,x_{\dell i})\}}, 
	   	\eqno(7.6)$$
where $\alpha_i(1),\ldots,\alpha_i(k)$ are specified or
unknown parameters, and
$\beta'H_i(x_i,x_{\dell i})
	 =\sum_{u=1}^p\beta_uH_{i,u}(x_i,\allowbreak x_{\dell i})$
for certain component functions 
$H_{i,u}$ that depend both on the 
$x_i$ at time position $i$ and of its neighbouring values $x_{\dell i}$. 
For a second order Markov chain, for example, one would typically have 
$H_i$ equal to a common $H$ function for $2\le i\le n-2$ and some
special functions at the borders. 
Maximum pseudo-likelihood estimation maximises 
${\rm PL}_n(\beta)=\prod_{i=0}^n f_\beta(x_i|x_{\dell i})$ 
w.r.t.~the parameters. 
See Hjort and Omre (1993, section 3.2), for example, 
for comments on this model building machinery in dimensions 1 and 2,
and for some comments on the difference between maximum PL and maximum 
likelihood. 

We may incorporate the $\alpha_i(x)$'s in the vector 
$H_i(x,x_{\dell i})$ of $H_{i,u}$-functions, for notational 
convenience. From Lemma A2 $\log {\rm PL}_n(\beta)$ is concave in $\beta$.
Consider 
$h_i(x_{\dell i})=\sum_{j=1}^kH_i(j,x_{\dell i})\,f_{\beta_0}(j|x_{\dell i})$
and
$$V_i(x_{\dell i})=\sum_{j=1}^k\bigl(H_i(j,x_{\dell i})-h_i(x_{\dell i})\bigr)
  	\bigl(H_i(j,x_{\dell i})-h_i(x_{\dell i})\bigr)'
  	\,f_{\beta_0}(j|x_{\dell i}), $$
which can be interpreted as respectively 
$\E\{H_i(X_i,x_{\dell i})|x_{\dell i}\}$ 
and ${\rm VAR}\,\{H_i(X_i,x_{\dell i})|x_{\dell i}\}$. 
After some work exploiting Lemma A2 one finds that 
$$\log{\rm PL}_n(\beta_0+s/\sqrt{n})-\log{\rm PL}_n(\beta_0)
=U_n's-\half s'(J_n/n)s-r_n(s), \eqno(7.7)$$
where 
$$U_n=n^{-1/2}\sum_{i=0}^n\{H_i(X_i,X_{\dell i})-h_i(X_{\dell i})\}
	\quad{\rm and}\quad
  J_n=\sum_{i=0}^nV_i(X_{\dell i}), $$	
and where in fact $r_n(s)=O(n^{-1/2})$. 
The usual arguments now give 
$\sqrt{n}(\hatt\beta_n-\beta_0)\arr_d\normal\{0,J^{-1}\}$
under mild assumptions, provided the assumed model (7.6) is correct.
Here $J$ turns out to be both the limit of $J_n/n$ as well as 
the covariance matrix in the limiting distribution for $U_n$. 
There is also an appropriate sandwich generalisation with
covariance matrix of type $J^{-1}KJ^{-1}$ outside model conditions. 
Doing the details here properly calls for a central limit theorem
and a weak law of large numbers for Markov chains,
and such can be found in Billingsley (1961), for example. 

These Markov random field models are more important 
in the 2- and 3-dimensional cases, where one enters the world 
of statistical image analysis. The method above can be used to prove 
consistency of the maximum PL estimator. 

%

\bigskip
{\bf Appendix.} 
Here we give three lemmas that were used at various stages above.
They should also have some independent interest.

\subsection
{\sl A1. Necessary and sufficient conditions for asymptotic
normality of linear combinations of binomials.}
The following result with further consequences was used in section 5.

{{\smallskip\sl
{\csc Lemma A1.} 
Consider independent Bernoulli variables $Y_{n,i}\sim{\rm Bin}\{1,q_{n,i}\}$,
and real numbers $z_{n,i}$ standardised to have 
$\sumin z_{n,i}^2q_{n,i}(1-q_{n,i})=1$. Then 
$\sumin z_{n,i}(Y_{n,i}-q_{n,i})\arr_d\normal \{0,1\}$
if and only if 
$$N_n(\delta)=\sumin z_{n,i}^2q_{n,i}(1-q_{n,i})\,
	I\{|z_{n,i}|\ge\delta\}\arr0
	\quad {\sl for\ each\ positive\ }\delta. \eqno(\hbox{\rm A.1})$$
\smallskip}}

{\csc Proof:} 
The Lindeberg condition is that 
$$\eqalign{
L_n(\delta)&=\sumin Ez_{n,i}^2(Y_{n,i}-q_{n,i})^2
	\,I\{|z_{n,i}(Y_{n,i}-q_{n,i})|\ge\delta\} \cr
	&=\sumin z_{n,i}^2q_{n,i}(1-q_{n,i})\,
	\bigl[q_{n,i}I\{|q_{n,i}z_{n,i}|\ge\delta\}
		+(1-q_{n,i})I\{|(1-q_{n,i})z_{n,i}|\ge\delta\}\bigr]}$$ 
should tend to zero for each positive $\delta$. 
It is not difficult to establish 
$\half N_n(2\delta)\le L_n(\delta)\le N_n(\delta)$,
so (A.1) is in fact equivalent to the Lindeberg requirement. 
In particular (A.1) implies a $\normal\{0,1\}$ limit.

Necessity is harder. Assume a $\normal\{0,1\}$ limit in distribution.
We first symmetrise in the following fashion: 
Let $\tilda Y_{n,i}=Y_{n,i}-Y_{n,i}'$ where 
$Y_{n,1}',Y_{n,2}',\ldots$ are independent
copies of $Y_{n,1},Y_{n,2},\ldots$, and let 
$$Z_n=\sumin z_{n,i}(Y_{n,i}-q_{n,i}),\quad 
	Z_n'=\sumin z_{n,i}(Y_{n,i}'-q_{n,i}),
	\quad {\rm and\ }\tilda Z_n=Z_n-Z_n'.$$
By assumption $\tilda Z_n\arr_d\normal\{0,2\}$. 
We first show that 
$$m_n=\maxin\,\min\{|z_{n,i}|,q_{n,i},1-q_{n,i}\}\arr0. $$
Otherwise there would be some $\eps>0$ such that say $|z_{n,1}|\ge\eps$
and $\eps\le q_{n,1}\le1-\eps$. 
Break $\tilda Z_n$ into a sum of $V_n=z_{n,1}\tilda Y_{n,1}$
and $W_n$, two independent and symmetric variables. 
Uniform tightness of $\tilda Z_n$ and
symmetry imply uniform tightness of both $V_n$ and $W_n$. 
Along some subsequence we would have $V_n\arr_dV$
and $W_n\arr_dW$, independent with $V+W$ distributed as $\normal\{0,2\}$.
By Cram\'er's theorem about convolution factors of $\normal\{0,2\}$ 
we would have $V$ normal. But $V$ is not degenerate, and cannot be
normal after all, since $V_n$ takes only three values. This proves
$m_n\arr0$.

But this implies the usual infinitesimal array property
$$\maxin{\rm Pr}\{|z_{n,i}\tilda Y_{n,i}|\ge\delta\}\arr0
	\quad {\rm for\ each\ }\delta. $$
For if $|z_{n,i}|<\delta$ then the probability is zero, and if 	
$|z_{n,i}|\ge\delta$ then the probability is 
$2q_{n,i}(1-q_{n,i})\le2m_n$
when $n$ is large enough for $m_n<\delta$ to hold.
Next look at page 92 of Petrov (1975). From limiting normality follows
$$\sumin{\rm Var}\,\big[z_{n,i}\tilda Y_{n,i}
	I\{|z_{n,i}\tilda Y_{n,i}|<\delta\}\big]\arr2. $$
If $|z_{n,i}|\ge\delta$ the indicator here picks out $\tilda Y_{n,i}=0$,
and there is no contribtion to the sum, whereas if $|z_{n,i}|<\delta$
the summand is $2z_{n,i}^2q_{n,i}(1-q_{n,i})$. Hence 
$\sumin z_{n,i}^2q_{n,i}(1-q_{n,i})\,I\{|z_{n,i}|<\delta\}\arr1$, 
and $N_n(\delta)\arr0$ follows from the assumed 
$\sumin z_{n,i}^2q_{n,i}(1-q_{n,i})=1$. \square 

\smallskip
The surprising thing here is that we do not need to explicitly 
assume $\maxin \E\{z_{n,i}(Y_{n,i}-q_{n,i})\}^2\arr0$, as with 
Feller's partial converse to the Lindeberg theorem;
it follows from asymptotic normality and the special properties
of the $Y_{n,i}$ sequence.  

Lemma A1 can next be used to address the vector case,
via the Cram\'er--Wold theorem. We phrase the result as follows,
to suit the development of section 5.  
If $x_1,x_2,\ldots$ is a sequence of $p$-vectors, 
and $Y_1,Y_2,\ldots$ are Bernoulli with $q_1,q_2,\ldots$, then 
$$J_n^{-1/2}\sumin (Y_i-q_i)x_i\arr_d\normal_p \{0,I_p\}, 
	\quad {\rm where\ }J_n=\sumin q_i(1-q_i)x_ix_i', 
		\eqno(\hbox{\rm A.2})$$
if and only if 
$$N_n(\delta)=\sumin x_i'J_n^{-1}x_i\,q_i(1-q_i)\,
	I\{|J_n^{-1/2}x_i|\ge \delta\}\arr0 
	\quad {\rm for\ each\ }\delta>0. \eqno(\hbox{\rm A.3})$$
This is proved by noting first that (A.2) is equivalent to 
$$N_n^0(s,\delta)=\sumin s'x_i'J_n^{-1}x_is\,q_i(1-q_i)\,
I\{|s'J_n^{-1/2}x_i|\ge\delta\} \arr0$$
for all $s$ with length 1 and all positive $delta$.
But $N_n^0(s,\delta)\le N_n(\delta)\le 
p^2\max_{j\le p}N_n^0(e_j,\delta/\sqrt{p})$, where $e_j$ 
is the $j$th unit vector.

A simple {\it sufficient} condition for (A.1) to hold is that 
$\lambda_n^0=\maxin |z_{n,i}|\arr0$, since 
the left hand side of (A.1) is bounded by 
$\lambda_n^0/\delta$. 
Similarly condition (A.3) is implied by the simpler condition 
$\lambda_n=\maxin |J_n^{-1/2}x_i|\arr0$, 
since $N_n(\delta)\le\lambda_np/\delta$.

\subsection
{\sl A2. Expansion lemma.}
The following result was used in section 6 and in 7A.

{{\smallskip\sl
{\csc Lemma A2.} 
(i) Suppose $K(t)=\log R(t)$, where $R(t)=\sumin w_i\exp(a_it)$
for certain nonnegative weights $w_i$, not all equal to zero, 
and arbitrary constants $a_i$. 
Let $v_i(t)=w_i\exp(a_it)/R(t)$ be the tilted and 
normalised weights, summing to one. 
Then $K(t)$ is convex with derivatives 
$$\eqalign{
K'(t)&=\sumin v_i(t)a_i=\bar a(t), \cr
K''(t)&=\sumin v_i(t)\{a_i-\bar a(t)\}^2, \cr
K'''(t)&=\sumin v_i(t)\{a_i-\bar a(t)\}^3. \cr}$$
(ii) The expansion  
$$\log\Bigl\{\sumin w_ie^{a_it}\Bigr\}-\log \Bigl\{\sumin w_i\Bigr\}
	=\bar a(0)t+\half\sumin v_i(0)\{a_i-\bar a(0)\}^2t^2
		+v(t) $$
holds, featuring untilted weights $v_i(0)=w_i/\sumin w_i$, 
with the following valid bounds on the remainder:
$$|v(t)|\le \hbox{$4\over3$}\mu_n^3|t|^3, \qquad
  |v(t)|\le \hbox{$2\over3$}
	g(\mu_n|t|)\sumin v_i(0)\{a_i-\bar a(0)\}^2|t|^2. 
		\eqno(\hbox{\rm A.4})$$
Here $\mu_n=\maxin |a_i-\bar a(0)|$ and $g$ is the function 
$g(u)=u\exp(2u+4u^2)$. 
\smallskip}}

{\csc Proof:} 
The formulae for the derivatives are proved 
by direct differentiation and inspection,
and convexity follows of course from the nonnegative second derivative. 
To prove (ii), consider the exact third order Taylor expansion 
$K(t)-K(0)=K'(0)t+\half K''(0)t^2+\sixth K'''(s)t^3$ 
for some suitable $s$ between $0$ and $t$. 
The problem is to bound the remainder term 
in terms of $\mu_n$. 

The first bound is easy. It follows upon observing that 
$|\bar a(s)-\bar a(0)|=|\sumin v_i(s)\{a_i-\bar a(0)\}|\le \mu_n$ 
and its triangle inequality consequence $|a_i-\bar a(s)|\le2\mu_n$,
since this yields $|K'''(s)|\le(2\mu_n)^3$. While this bound 
often suffices we shall have occasion to need the sharper 
second bound too. The point is to exploit the fact that 
$s$ is bounded by $|t|$ when bounding $|K'''(s)|$. 
Start out writing $v_i(s)=v_i(0)(1+\eps_i)$, where 
some analysis shows that $\exp(-2\mu_n|t|)\le 1+\eps_i\le \exp(2\mu_n|t|)$.
Then 
$$\eqalign{
|\bar a(s)-\bar a(0)|
&=|\sumin v_i(0)(1+\eps_i)\{a_i-\bar a(0)\}| \cr
&=|\sumin v_i(0)\eps_i\{a_i-\bar a(0)\}| 
 \le K''(0)^{1/2}\Bigl\{\sumin v_i(0)\eps_i^2\Bigr\}^{1/2}. \cr}$$
A further bound on the right hand side is 
$K''(0)^{1/2}\delta_n=K''(0)^{1/2}\maxin |\eps_i|$. This gives 
$$\eqalign{
|K'''(s)|&\le 2\mu_n\sumin v_i(s)|a_i-\bar a(s)|^2 \cr
	 &\le 4\mu_n\sumin v_i(0)(1+\eps_i)
		\{|a_i-\bar a(0)|^2+|\bar a(0)-\bar a(s)|^2\} \cr
	 &\le 4\mu_n(1+\delta_n)(1+\delta_n^2)K''(0). \cr}$$
But some checking reveals 
$1+\delta_n\le\exp(2\mu_n|t|)$ and $1+\delta_n^2\le\exp(4\mu_n^2|t|^2)$. 
This shows $|K'''(s)|\le4|t|^{-1}g(\mu_n|t|)K''(0)$ with 
the $g$-function given above. \square 

\subsection
{\sl A3. Dominated convergence theorem for convergence in probability.}
This result was used several times in section 6, and a close relative
was used in 4B. 

{{\smallskip\sl
{\csc Lemma A3.}
Let $0\le X_n(s,\omega)\le Y_n(s,\omega)$ be jointly $(s,\omega)$-measurable 
random functions on the interval $[0,L]$.
Suppose $\lambda$ is a measure such that
$Y_n(s)\arr_p Y(s)$ and $X_n(s)\arr_p X(s)$ for $\lambda$ almost all $s$
and that $\int Y_n(s)\,\d\lambda(s)\arr_p\int Y(s)\,\d\lambda(s)$,
a limit finite almost everywhere.
Then $\int X_n(s)\,\d\lambda(s)\arr_p\int X(s)\,\d\lambda(s)$ too.
\smallskip}} 

{\csc Proof:} 
It is enough to check almost sure convergence for a subsubsequence of each
subsequence. By convergence in probability (for $\omega$, with $s$ fixed) 
and then dominated convergence, we have
$$\pi_n(\epsilon) :=
	(\PP\otimes\lambda)\{(\omega,s)\colon
	|X_n(\omega,s)-X(\omega,s)|> \eps\} \to 0. $$
A similar result holds for $\{Y_n\}$.  
Replace $\eps$ by a sequence
$\{\eps_n\}$ decreasing to zero, then extract a subsubsequence
along which the sequence of integrals is convergent.
For some set $N$ with $(\PP\otimes\lambda)N=0$ we get convergence for all
$(\omega,s)\in N^c$ of both the $X_n$ and the $Y_n$ subsubsequences.  
For almost all $\omega$, therefore,  
$\lambda\{s\colon(\omega,s)\in N\}=0$. 
Finally argue using the Fatou Lemma for $Y_n\pm X_n$ along the 
subsubsequences to get the result. \square 


\smallskip
This is nice in that it circumvents the need to establish
uniformity of the convergence in probability; 
this is typically more difficult to ascertain than pointwise
convergence in probability. 
The lemma was used several times in sections 6 and 7A, 
partly in the form of the following useful corollary: 
If in particular $Z_n(s)\arr_p 0$ for almost all $s$, then 
$\int_0^L Y_n(s)\,I\{|Z_n(s)|\ge\delta\}\,\d{}s\arr_p 0$.
It can also be used to simplify the Lindeberg type condition 
in the form of Rebolledo's martingale convergence theorem 
given in Andersen and Gill (1982, appendix). 

\bigskip
{\bf Acknowledgements.}
This paper resulted from a collaboration made possible by the support
of the Mathematical Sciences Research Institute at Berkeley during the
Autumn of 1991, under National Science Foundation grant 8505550. We are
grateful to the MSRI for the opportunity to learn from each other.
N.L.H.~has also been partly supported by grants 
>from the Royal Norwegian Research Council
and D.P.~partly supported by National Science Foundation Grant DMS-9102286.

\bigskip\bigskip
\centerline{\bf References}

\parindent0pt
\baselineskip11pt
\parskip3pt 
  
\def\ref#1{{\noindent\hangafter=1\hangindent=20pt
  #1\smallskip}}  

\medskip
\ref{%
Andersen, P.K., Borgan, \O., Gill, R.D., and Keiding, N. (1992).
{\sl Statistical Models Based on Counting Processes.}
Springer-Verlag, New York.} 

\ref{%
Andersen, P.K.~and Gill, R.D. (1982).
Cox's regression model for counting processes:
a large sample study. 
{\sl Ann.~Statist.}~{\bf 10}, 1100--1120.} 

\ref{%
Bassett, G.W.~and Koenker, R.W. (1982).
An empirical quantile method for linear models with iid errors.
{\sl J.~Amer.~Statist.~Assoc.}~{\bf 77}, 407--415.}   

\ref{%
Bickel, P.J., Klaassen, C.A., Ritov, Y., and Wellner, J.A. (1992).
{\sl Efficient and Adaptive Inference in Semiparametric Models.}
To exist.}

\ref{%
Billingsley, P. (1961). 
{\sl Statistical Inference for Markov Processes.}
University of Chicago Press, Chicago.}

\ref{%
Borgan, \O. (1984).
Maximum likelihood estimation in parametric counting process models,
with applications to censored failure time data. 
{\sl Scand.~J.~Statist.}~{\bf 11}, 1--16.
[Corrigendum, ibid.~p.~275.] }




\ref{%
Haberman, S.J. (1989).
Concavity and estimation.
{\sl Ann.~Statist.}~{\bf 17}, 1631--1661. }

\ref{%
Hjort, N.L. (1986). 
Bayes estimators and aymptotic efficiency in parametric counting
process models. 
{\sl Scand.~J.~Statist.}~{\bf 13}, 63--85. }


\ref{%
Hjort, N.L. (1988a). 
Logistic regression when the model is wrong. 
Appendix to {\sl Statistical models for the probability
of finding oil or gas}, with V.~Berteig. 
Norwegian Computing Centre Report, Oslo.}


\ref{%
Hjort, N.L. (1988b).
Lecture notes and seventy-five exercises
on resampling and bootstrapping.
Statistical Research Report, University of Oslo.}

\ref{%
Hjort, N.L. (1991).
Semiparametric estimation of parametric hazard rates.
In {\sl Survival Analysis: State of the Art}, 
Kluwer, Dordrecht, pp.~211--236. 
Proceedings of the {\sl NATO Advanced Study Workshop on Survival
Analysis and Related Topics}, Columbus, Ohio, eds.~P.S. Goel and J.P.~Klein.} 

\ref{%
Hjort, N.L. (1992).
On inference in parametric survival data models.
{\sl Intern.~Statist.~Review} {\bf 60}, 355--387. }

\ref{%
Hjort, N.L.~and Omre, H. (1993).
Topics in spatial statistics (with discussion). 
{\sl Scand.~J.~Statist.}, to appear.  }

\ref{%
Huber, P. (1967).
The behaviour of maximum likelihood estimates under non-standard conditions.
{\sl Proceedings of the Fifth Berkeley Symposium 
on Mathematical Statistics and Probability} {\bf 1}, 221--233.
University of California Press. }



\ref{%
Jure\v ckov\'a, J. (1977). 
Asymptotic relations of M-estimates and R-estimates in linear regression model.
{\sl Ann.~Statist.}~{\bf 5},464--472.}

\ref{%
Jure\v ckov\'a, J. (1992). 
Estimation in a linear model based on regression rank scores.
{\sl J.~Nonparametric Statistics} {\bf 1}, 197--203.}

\ref{%
Lehmann, E.L. (1983). 
{\sl Theory of Point Estimation.}
Wiley, New York.}

\ref{%
McCullagh, P.~and Nelder, J.A. (1989).
{\sl Generalized Linear Models} (2nd edition).
Chapman \& Hall, London. }

\ref{%
Niemiro, W. (1992).
Asymptotics for M-estimators defined by convex minimization.
{\sl Ann.~Statist.} {\bf 20}, 1514--1533. }


\ref{%
Petrov, V.V. (1975).
{\sl Sums of Independent Random Variables.}
Springer-Verlag, New York. }

\ref{%
Pollard, D. (1984).
{\sl Convergence of Stochastic Processes.}
Springer-Verlag, New York.}

\ref{%
Pollard, D. (1985).
New ways to prove central limit theorems.
{\sl Econometric Theory} {\bf 1}, 295--314.}

\ref{%
Pollard, D. (1990).
{\sl Empirical Processes: Theory and Applications.}
NSF-CBMS Regional Conference Serires in 
Probability and Statistics, Vol.~2.
IMS, Hayward, California. }

\ref{%
Pollard, D. (1991).
Asymptotics for least absolute deviation regression estimators.
{\sl Econometric Theory} {\bf 7}, 186--199. }

\ref{%
Rockafellar, R.T. (1970).
{\sl Convex Analysis.} 
Princeton University Press, Princeton, New Jersey.}


\bye